\documentclass
[final]{siamltex}
\usepackage{titlesec}

\titleformat{\paragraph}
{\normalfont\normalsize\bfseries}{\theparagraph}{1em}{}
\titlespacing*{\paragraph}
{0pt}{3.25ex plus 1ex minus .2ex}{1.5ex plus .2ex}
\newcommand{\Z}{\mbox{\rm \lower0.3pt\hbox{$\angle\!\!\!$}Z}}

\newtheorem{remark}{Remark}[section]
\newcommand{\li}{\infty}
\newcommand{\xnp}{x_{n+1} }
\newcommand{\xn}{x_{n} }

\newcommand{\lb}{\left (}
\newcommand{\rb}{\right )}
\newcommand{\la}{\langle }
\newcommand{\ra}{\rangle }

\newcommand{\ddem}{\ind  {\it Proof.} }

\newcommand{\edem}{. $\Box$}

\newcommand{\ind}{\hspace*{0.3cm} }

\newcommand{\bc}{\begin{center}}
\newcommand{\ec}{\end{center}}
\newcommand{\benq}{\begin{eqnarray}}
\newcommand{\be}{\begin{equation}}

\newcommand{\ee}{\end{equation}}
\newcommand{\eenq}{\end{eqnarray}}
\newcommand{\ba}{\begin{array}{l}}
\newcommand{\ea}{\end{array}}

\newcommand{\bfg}{\begin{figure}[h]}
\newcommand{\efg}{\end{figure}}

\newcommand{\un}{\underline}
\newcommand{\mul}{\times}
\newcommand{\esc}{\hskip 0.1cm}
\newcommand{\se}{\hskip 0.1cm}
\newcommand{\es}{\hskip 0.3cm}

\newcommand{\al}{\alpha}
\newcommand{\ld}{\lambda}
\newcommand{\ldn}{\lambda_n}

\newcommand{\Hs}{{\cal H}}

\newcommand{\ynp}{y_{n+1} }
\newcommand{\yn}{y_{n} }
\newcommand{\dxnp}{\dot{x}_{n+1} }
\newcommand{\dxn}{\dot{x}_{n} }

\newcommand{\dzn}{\dot{z}_{n} }

\newcommand{\dynp}{\dot{y}_{n+1} }
\newcommand{\dyn}{\dot{y}_{n} }

\newcommand{\xnm}{x_{n-1} }
\newcommand{\vn}{v_{n} }
\newcommand{\vnp}{v_{n+1} }
\newcommand{\dvnp}{\dot{v}_{n+1} }
\newcommand{\ttn}{\theta_n }

\newcommand{\der}{\partial}

\newcommand{\gamos}{\vartheta_{n}}
\newcommand{\gamc}{\xi}
\newcommand{\bgamc}{(1-\gamc_n)}

\usepackage{amssymb}
\title{FAST CONVERGENCE OF  GENERALIZED  FORWARD-BACKWARD ALGORITHMS  FOR STRUCTURED MONOTONE INCLUSIONS}
\author{
 Paul-Emile Maing\'e\thanks{Universit\'e des Antilles,
              D.S.I.,  Campus de Schoelcher, 97233 Cedex, Martinique,  F.W.I.,  MEMIAD
   ({\tt Paul-Emile.Mainge@univ-antilles.fr})}
   }
\usepackage{graphics}
\usepackage{latexsym}
\usepackage{subeqn}
\usepackage{cite}
\input epsf
\date{555}
\begin{document}
\maketitle
%
\begin{abstract}
 In this paper, we develop rapidly convergent    forward-backward   algorithms for computing  zeroes  of the sum
  of  finitely many maximally monotone operators. A modification of the classical  forward-backward method for two general operators is first considered, by   incorporating   an inertial term   (closed to the acceleration techniques introduced by Nesterov), a constant relaxation factor  and  a correction  term. In a Hilbert space setting, we prove  the weak  convergence to equilibria of the iterates  $(x_n)$, with     worst-case  rates of $ o(n^{-2})$ in terms of both the  discrete velocity  and the  fixed point residual, instead of the   classical rates of ${\cal O}(n^{-1})$  established so far for related algorithms. Our procedure  is then adapted  to more  general monotone inclusions and a fast primal-dual  algorithm  is proposed for solving convex-concave saddle point problems.
\end{abstract}



\pagestyle{myheadings} \thispagestyle{plain} \markboth{P.E.
MAINGE}{FAST CONVERGENCE OF GENERALIZED  FORWARD-BACKWARD ALGORITHMS}
\section{\large Introduction.}
\setcounter{equation}{0}
Let    ${\cal H} $   be  a  real Hilbert space endowed
with inner product
  and   induced norm denoted by  $\la .,. \ra$ and  $\|\se.\se\|$, respectively. For any given  linear self-adjoint and positive definite  mapping ${\cal M }: \Hs \to \Hs$, and any bounded linear operator ${\cal L }: \Hs \to \Hs$,  we set
  $ \|\se.\se\|_{\cal M}= \sqrt{\la {\cal M} (.),. \ra}$ (as an auxiliary metric on $\Hs$) and  $\|{\cal L}\|= \sup_{\|x\|=1} \|{\cal L} x\|$.
    Our goal  is to propose and study    rapidly converging  forward-backward  methods  for solving a wide class of structured monotone inclusions, but it can be applied to a much larger  class of problems.
 \subsection{A new way to speed up  forward-backward methods.}
   The simplest   case of the considered  problems  consists of finding a zero of the sum of two   operators, especially the following  monotone inclusion
 \be \label{pbg} \mbox{find $\bar{x} \in \Hs$ such that $0 \in A (\bar{x})  +B (\bar{x}) $}, \ee
 under the  assumptions that $A: \Hs \to 2^{\Hs}$ is maximally monotone and that $B : \Hs \to {\Hs}$  is co-coercive (see, e.g.,  \cite{baco, combvu}).
  This framework  finds many important applications in scientific fields such as image processing, computer vision, machine learning, signal processing, optimization, equilibrium theory,
economics and game theory, partial differential equations, statistics, among other subjects
(see, e.g., \cite{baco, comb, lp,lmh,tak, attsou, rocki,zb}).
In particular, we  recall that the above setting  encompasses  the non-smooth structured convex minimization problem
 \be
 \label{pbc}
 \min_{\Hs}  \{\Theta:=f+g\},
 \ee
 where $g: \Hs \to (-\li, \li]$  is  proper,  convex and lower semi-continuous and $f: \Hs \to (-\li, \li)$  is   convex  differentiable with  Lipschitz continuous gradient.
  Indeed,  the gradient of a convex and Frechet differentiable function is well-known to be  co-coercive provided that it is Lipschitz continuous (see \cite{baco}). Thus,  (\ref{pbc})  is nothing but  the special instance  of  (\ref{pbg}) when $A=\der f$ ($\der f$ being   the Fenchel sub-differential of $f$) and $B=\nabla g$ ($\nabla g$ being    the gradient of $g$). \\
\ind  A typical  method  for solving  (\ref{pbg}), with a $\beta$-co-coercive operator $B$ (for some $\beta>0$),  is the forward-backward algorithm (FBA, for short), which operates according to the routine (see, e.g.,  Lions-Mercier \cite{lm}, Passty \cite{pass})
\be \label{fbpa} \xnp= J_{\ld_n  A}(\xn- \ldn  B(\xn)), \ee
 where $ (\ld_n) \subset    (0,2\beta)$ and $J_{\ld_n  A}:=(I+\ld_n A)^{-1}$ is the resolvent operator of $A$ with index $\ld_n$ (see Section \ref{resol} for more details on the resolvent operators). \\
 \ind  When applied to the minimization problem (\ref{pbc}), FBA  was shown to generate weakly convergent sequences $(x_n)$, with worst-case rates $\Theta(x_n)-\inf_{\Hs} \Theta ={\cal O}(n^{-1})$ (for the function values)  and    $\|\xnp-\xn\|^2={\cal O}(n^{-1})$ (for the discrete velocity). These rates   have  been considerably enhanced  through an  accelerated
forward-backward algorithm  (AFBA, for short, which will be described in   Section \ref{smap}),   by means of an inertial-type extrapolation process. The latter AFBA   generates  convergent sequences $(x_n)$  with  the improved worst-case  rates  $\Theta(x_n)-\inf_{\Hs} \Theta =o(n^{-2})$ and
 $\|\xnp-\xn\|^2= o(n^{-2})$. Afterwards,  general variants  of AFBA,  with  arbitrary  monotone  operators $A: \Hs \to 2^{\Hs}$ and $B: \Hs \to \Hs$,  have  been extensively studied and many attempts to get convergence rates  similar to that  reached for the convex minimization  case   can be found in the literature (see, e.g., \cite{ac, lp, moudol}). However,   most of these works  are only concerned with  empirical results.  \\
 \\
 \ind Our purpose here is twofold with regard to the previous observations:\\
  \\
  \ind First, we aim at  extending,  in terms of discrete velocity and fixed point residual,   the above convergence properties (obtained for  structured  convex minimization)  to some  new variant  of (\ref{fbpa}). This process  will be investigated with an additional  preconditioning strategy, as  suggested  by  Lorenz-Pock \cite{lp},   under the following general  conditions:
    \begin{subeqnarray}
    \label{cdi}
     \label{cdi1}
    && \mbox{$A : \Hs \to 2^{\Hs}$ is   maximally monotone on $\Hs$};  \\
     \label{cdi2}
     && \mbox{ $B : \Hs \to {\Hs}$ is  co-coercive w.r.t a linear, self-adjoint and positive definite map  $L$}, \\
    \label{cdi3}
    && \mbox{  $S:=(A+ B)^{-1}(0) \neq \emptyset$}.
     \end{subeqnarray}
     %
    \ind An operator   $B : \Hs \to 2^{\Hs}$ is said to be  co-coercive w.r.t a linear and positive definite   map  ${L}: \Hs \to \Hs $ if it satisfies
     \be \label{dcc} \mbox{ $ \la Bx-By, x-y \ra \ge  \|Bx-By\|_{L^{-1}}^2$, \es for $(x, y) \in \Hs^2$}. \ee
    \ind Note that in the simple case when $L=\beta^{-1}I$ (for some  $\beta >0$), any  operator $B$ verifying (\ref{dcc}) is $\beta$-co-coercive. It can also turn out useful   to consider more general map $L$ as discussed in \cite{lp}. \\
    \\
\ind Specifically, for solving (\ref{pbg})-(\ref{cdi}), we introduce  CRIFBA (corrected relaxed   inertial forward-backward algorithm) which consists   of
     sequences $\{z_n, x_n\} \subset \Hs$ generated  by the following process:\\
\\
 {\bf (CRIFBA)}: \\
 $\rhd$ {\bf Step 1} (\un{initialization}): \\
 \ind \ind Let $M: \Hs \to \Hs$ be  a linear self-adjoint and positive  definite map, \\
 \ind \ind let $\{ z_{-1}, x_{-1}, x_0 \} \subset \Hs$,  $\{s_1, s_0, \nu_0 \} \subset [0,\li)$,  $\{e,  \ld, w \} \subset (0,\li)$, \\   \ind \ind and set
 \be  \label{vvtn}
  \mbox{$\nu_n=s_1n+\nu_0$, \es $\ttn= 1- \frac{e+s_1 }{e+ \nu_{n+1}}$, \es $\gamma_n=  1-\frac{s_0}{e+ \nu_{n+1}}$}. \ee
 $\rhd$ {\bf Step 2} (\un{main step}): \\
\ind \ind Given $\{z_{n-1},x_{n-1}, x_n\} \subset \Hs$ (with $n \ge 0$), we   compute the updates by
 \begin{subeqnarray}
   \label{sc}
\label{sc1a}
&  & \hskip -3.cm \mbox{$v_n= z_{n-1}- \xn$},\\
\label{sc1c}
&  & \hskip -3.cm \mbox{$z_n= x_n + \ttn   (\xn -\xnm) +  \gamma_n  v_n $},\\
\label{sc1b}
& &  \hskip -3.cm \mbox{$\xnp=(1-w) z_n + w J_{\ld M^{-1} A} (z_n -\ld M^{-1} B(z_n)) $}.
  \end{subeqnarray}     
  where $J_{\ld  M^{-1} A}:=(I+\ld  M^{-1} A)^{-1}$ denotes the resolvent  of $ M^{-1}A$  (also referred to as a generalized resolvent of $A$).\\
  \\
   \ind The  algorithm under consideration can be regarded as a preconditioned and relaxed variant of  the classical forward-backward    in which we further incorporate the momentum term "$\ttn (\xn-\xnm)$"  (inspired by  Guler's acceleration techniques for convex minimization) and  the    correction term "$\gamma_n (z_{n-1}-x_n)$" (similar to that suggested by Kim \cite{kim} in proximal point iterations).   Note that (\ref{sc1b}) can be reformulated as
  \be \xnp=(1-w) z_n + w (M+ \ld  A)^{-1}(M z_{n}-\ld  B(z_{n})). \ee
 \ind Then, as illustrated in  \cite{lp} in the context of convex-concave saddle-point problems, the introduction of the map $M$   can be helpful in some situations  to make the iteration feasible, but it can be also interpreted as a left pre-conditioner to the monotone inclusion (\ref{pbg}). Furthermore, this procedure   allows us to investigate our algorithm with respect to the auxiliary metric $\|.\|_M$, which can be also used as a speeding up process.\\
  \\
 \ind The simple form of CRIFBA allows us to extend   the fast   convergence rates obtained  for  AFBA (in the context of potential operators)  to the   wide framework of the structured monotone inclusion (\ref{pbg})-(\ref{cdi}),  especially in terms of the discrete velocity and the fixed point residual  $\|G^M_{\ld}(x_n)\|_M^2$, where
  \be \label{defmld} \mbox{$G^{M}_{\ld}= \frac{1}{\ld } \lb I-J_{\ld M^{-1} A}\circ (I-\ld M^{-1} B)  \rb $}, \ee
 hence $\|G^M_{\ld}(x_n)\|_M=(\ld w)^{-2}\|x_n-z_{n-1}\|_M^2$.  The term   $G_{\ld}^M(z)$ (for $z \in \Hs$) can be regarded as a tool to measure the  accuracy of some point $z$  to the solution set  $S$ since $ G_{\ld} (z)=0$ is equivalent to $z \in S$.
  In particular,  using appropriate  parameters $\{\theta_n, \gamma_n, \ld, w \}$,  we establish, among others (see, Theorem \ref{thmg}), the  weak convergence of $(x_n)$ towards equilibria belonging to $S$,  with the worst-case  rates
 $\|\xnp-\xn\|_M^2= o(n^{-2})$ (for the discrete velocity) and $\|G^M_{\ld}(x_n)\|_M^2=o(n^{-2})$ (for the  fixed point residual).
\\
  \\
  \ind Secondly, by  following the methodology developed by Raguet-Fadili-Peyré \cite{rafape}, we aim at  adapting  the proposed  acceleration techniques to the     more general
inclusion problem :
\be \label{gmoni} \mbox{find $\bar{x} \in \Hs$ s.t. $ 0 \in B(\bar{x}) +  \sum_{k=1}^p A_k(\bar{x})$,} \ee
where $B: \Hs \to \Hs$ is $\beta$-co-coercive  and  $(A_i)_{i=1}^p: \Hs \to 2^{\Hs}$ is a family of $p$  maximally  monotone operators. This gives rise to G-CRIFBA (generalized corrected and regularized forward-backward algorithm),
 whose structure  bears similarities with FBA. Indeed,  G-CRIFBA
 consists of an explicit forward step, followed by an implicit step in which the resolvent  operators of each   $A_i$ are computed in parallel. Let us underline that G-CRIFBA  provides yet another way for computing the resolvent of the sum of maximally
monotone operators at a point $y \in {\rm rang} \lb I+ \sum _{i=1}^p A_i \rb $. This can be seen  when  taking $B$ as the particular operator defined for $x \in \Hs$ by  $B(x) = x-y$ (hence $B$ is $1$-co-coercive).  A particular attention will also be paid to the case of potential operators. \\
\subsection{Motivations and reminders on proximal splitting algorithms.}
\subsubsection{\label{resol} Some splitting variants of the basic proximal  algorithm.}
%
 A classical method for
 computing zeroes of a maximally monotone operator  $A:\Hs \to 2^{\Hs}$ is the so-called  PPA  (proximal point algorithm) (see Martinet \cite{mart},
 Rockafellar \cite{roc, Ro-We}), which consists of  iteration
 \be \label{ppa} \xnp= J_{\ld  A }( \xn), \ee
 where $ J_{\ld  A }:=(I+\ld A)^{-1}$  (for some positive parameter $\ld$) is the resolvent operator of $A$, which is well-known to be everywhere defined and single-valued (see,   e.g., \cite{brez,gu,lem} for more details).
   In particular, when $A=\der g$ is the Fenchel sub-differential of a proper  convex and lower-semi-continuous function $g: \Hs \to(-\li,\li]$, the resolvent operator of $A$ reduces to the proximal mapping of $g$ given by
 \be \mbox{${\rm prox}_{\ld  g}(x):=J_{\ld \der g }={\rm argmin} _{ y \in \Hs} \lb g(y)+ (2\ld )^{-1}\|x-y\|^2 \rb$}. \ee
  It is well-known that  PPA   generates weakly convergent sequences $(x_n)$ with  a discrete velocity  that vanishes at the rate $\|\xnp-\xn\|^2= {\cal O}(n^{-1})$ (or equivalently,  $\|\xn - J_{\ld  A }(\xn)\|^2= {\cal O}(n^{-1})$ for the fixed point residual). \\
  \\
 %
 %
 %
 %
 \ind In many situations,  however, evaluating   the resolvent of the sum of two  maximally monotone    operators $A$ and $B$ turns out to be more complicated that evaluating separately the proximal operators of $A$ and $B$.    This observation  gave rise to two main categories  of  splitting methods of practical   interest (for instance, in the context of
sparse signal recovery \cite{daubdd,combwaj}, machine learning \cite{dusing} and image processing \cite{rafape}): \\
\\
\ind  {\rm \bf (c1)}  The first category   of splitting algorithms is composed of those that essentially include  backward  steps, namely,  no evaluations  of $A$ and  $B$, but only  evaluations of both  the resolvent  operators of  $A$ and $B$.  This kind of methods originates from   the Peaceman-Rachford and  the Douglas-Rachford splitting algorithms  (see \cite{dourach, peacerach, lionmerc}).  As an example, we mention the following iteration
(see  Corman-Yuan \cite{coryu}, Eckstein-Bertsekas \cite{ecbe})
\be \label{rmp} \mbox{$\xnp =    H (\xn) $, \ind where $H= J_{\ld A} ( 2  J_{\ld B}-I)+ I- J_{\ld B}$ \es (for some  $\ld >0$)}. \ee
  The operator defined by  $D:=H^{-1} -I$    was shown to verify   the remarkable properties  (see \cite{ecbe})
  \be
   \mbox{$D$  is maximally monotone, \ind
    $H= J_{D}$, \ind
    $D^{-1}(0)=S:=(A+B)^{-1}(0)$}.
   \ee
   Thus,   (\ref{rmp}) is nothing but the fixed point iteration for the resolvent  operator  $J_{D}$. It is a classical matter to see that the   iterates  $(x_n)$ produced by (\ref{rmp})     converge weakly to some element of $S$. \\
\\
 \ind  {\rm \bf (c2)}   The splitting methods in the  second category    combine both backward steps (evaluations of resolvent operators)  and forward steps (evaluations of one of the two operators). \\
 \ind - A popular example of such methods, for a single-valued mapping  $B$,  is given by so-called forward-backward algorithm (\ref{fbpa}). In specific, when $B$ is $\beta $- co-coercive,  (\ref{fbpa}) was shown to be  a convergent method, provided that    $(\ld_n)   \subset   (0, 2 \beta)$. Note indeed that this latter algorithm
 be reformulated  as  the fixed point iteration $\xnp = T_{\ld_n } x_n$, where
 $T_{\ld_n}=J_{\ld_n  A}\circ (I-\ld_n  B)$, with a fixed point set ${\rm Fix}( T_{\ld_n})=S:=(A+B)^{-1}(0)$.
  One can easily check that this assumption ensures  that  $T_{\ld_n}$ is non-expansive (since $ J_{\ld_n A}$ and  $I-\ld_n B$ are non-expansive).
   This classically leads to the weak convergence of
    the  iterates $(x_n)$ generated  by   (\ref{fbpa}) towards   some element of $S$.  \\
 \ind - Another example, for a single-valued mapping  $B$,  is given by the following  forward-backward-forward algorithm proposed  by
Tseng \cite{tseng} (see also \cite{botsedvu})
\begin{eqnarray}
\label{fbf}
  &  &  \mbox{$y_n=J_{\ld_n  A}(x_n-\ld_n B(x_n)), \es \es    \xnp =y_n -  \ld_n  (  B  (y_n)- B(x_n))$}.
 \end{eqnarray}
This  method   was shown to generate (weakly) convergent sequences, even when   $B$ is $L$-Lipschitz continuous (which is a weaker condition than  co-coerciveness), provided that    $(\ld_n)   \subset  (0,  L^{-1})$.\\
\\
\ind Later on,  generalized variants of the forward-backward  and forward-backward-forward algorithms (see, e.eg., \cite{combpes,rafape})  have   been adapted  more general inclusion problems such as  (\ref{gmoni}).\\
 \subsubsection{\label{smap} Proximal splitting  algorithms and acceleration processes.}
  It is suitable  to accelerate the proximal point method and its splitting variants, in view of  their various applications.   \\
  \\
\ind Note that,  in the context  of the convex minimization  (\ref{pbc}), the forward-backward method (\ref{fbpa})
 reduces to
\be \label{fba} \xnp= {\rm prox}_{\ld_n  g}(x_n-\ld_n \nabla f(x_n)). \ee
 For values $(\ld_n) \subset (0,L^{-1})$,  when $\nabla f$ is assumed to be $L$-Lipschitz continuous,  (\ref{fba})   generates (weakly) convergent sequences $(x_n)$ that satisfy the
(sub-linear) rate   $\Theta(x_n)-\inf_{\Hs} \Theta= {\cal O}(n^{-1})$.
   (\cite{bl}). \\
  \ind  Afterwards, (\ref{fba}) was enhanced  through the
 {\it Fast Iterative Soft Thresholding Algorithm} (FISTA) proposed by Beck-Teboule \cite{bteb} (also see \cite{visalbalver}), based  upon the acceleration techniques of  Guler \cite{gu2} and
 Nesterov \cite{nest, nest1, nest2}). FISTA was shown to produce iterates $(x_n)$ that  guarantee   a rate of convergence  $\Theta(x_n)-\inf_{\Hs} \Theta= {\cal O}(n^{-2})$. However,  the  convergence of these iterates   has  not been  established.  \\
 \ind This drawback was overcame by the following variant  of FISTA recently   introduced by  Chambolle-Dossal (see \cite{chamdos})
(also see  Attouch-Peypouquet \cite{attpey1})
   given  by
   \be \ba
  \label{afba}
  \label{afba1}
   z_{n} =x_n + \frac{n-1}{ n+\al-1} (\xn - x_{n-1}), \\
   \xnp= {\rm prox}_{\ld  g}( z_{n} -\ld  \nabla f(  z_{n} )),
   \ea \ee
  where  $\ld  \in  (0,  L^{-1})$ and $\al >0$.  It was proved for    $\al >3$  (see \cite{attpey1}) that    (\ref{afba})   generates (weakly) convergent sequences $(x_n)$ that minimize the function values $\Theta(x_n)$ with  a complexity result of $o(n^{-2})$, instead of the rates ${\cal O}(n^{-1})$ and ${\cal O}(n^{-2})$  obtained  for  ({\ref{fba}) and FISTA, respectively. \\
  \\
  \ind  It is worthwhile noticing in the case of  an arbitrary maximally monotone operator $A: \Hs \to 2^{\Hs}$ that  accelerated variants of PPA have  been  proposed and investigated   through RIPA (Regularized Inertial Proximal Algorithm)
by  Attouch-Peypouquet \cite{attpey},
PRINAM  (Proximal Regularized Inertial Newton Algorithm)  by  Attouch-Laslo \cite{attlas}. These algorithms, despite their interesting asymptotic features, require  unbounded proximal indexes for convergence and so cannot be extended to the forward-backward framework.
In the same context, an  accelerated proximal
point method involving constant proximal indexes was  proposed by  Kim \cite{kim}, based on the performance estimation problem (PEP) approach of Drori-Teboulle \cite{drori}.
   This yields the worst-case convergence rate of ${\cal O}(n^{-2})$  in terms of  fixed point residuals.
     Once again,  no convergence of the iterates was   established.
 To the best of our knowledge, regarding the existing algorithmic solutions to   (\ref{pbg})-(\ref{cdi}) with general operators, there are no analogous theoretical convergence results to that obtained for (\ref{afba}). Only somewhat empirical accelerations  have been proposed, except for the work of  Attouch-Cabot \cite{ac},  via  relaxation and  inertial extrapolation techniques. Some of these processes are recalled below:  \\
\ind {\rm \bf (e1)} Inertial  variants of (\ref{fbpa}), with  a co-coercive operator $B$,   have  been discussed by
Moudafi-Oliny \cite{moudol}:
\be
\label{moc}
 z_n= x_n + \al_n (\xn -\xnm), \ind
  \xnp= J_{\ld _n A} ( z_n - \ld _n B (x_n)),
\ee
  and by    Lorenz-Pock \cite{lp}:
\be
\label{poc}
 z_n= x_n + \al_n (\xn -\xnm), \ind
  \xnp= J_{\ld _n A} ( z_n - \ld _n B(z_n)),
\ee
where   $(\al_n)$ and ($\ld _n)$ are positive sequences. Note that this second method involves the evaluation of $B$ at $z_n$ instead of of $x_n$ (as done in  (\ref{moc})).\\
\ind  {\rm \bf (e2)}
A reflected variant of (\ref{fbpa}), with a Lipschitz continuous operator $B$,  was investigated  by  Cevher-Vu \cite{ce}:
 \be y_n= 2 \xn-\xnm , \es \xnp =  J_{\ld  A}(\xn-\ld  B\yn). \ee
\ind  {\rm \bf (e3)}  An  inertial and relaxed variant of  (\ref{fbpa}), with  a $\beta$-co-coercive operator $B$,  has  been discussed by Attouch-Cabot \cite{ac}:
\begin{eqnarray}
\label{pac}
 z_n= x_n + \al_n (\xn -\xnm), \es
 \xnp= (1-w_n) z_n + w_n J_{\ld_n A} ( z_n - \ld_n B z_n),
\end{eqnarray}
where $\{\al_n, w_n, \ld_n \}$ are positive and bounded sequences. Under various conditions on the  parameters, the authors have established   the weak convergence to equilibria of the iterates  $(\xn)$, with additional convergence rates   in terms of the discrete velocity and the  fixed point residual  $\|G^I_{\ld}(x_n)\|^2$   (where $G^I_{\ld}$ was  introduced in (\ref{defmld})). In particular,  this  work  meets the   setting  of Nesterov's accelerated methods  through the  choice
 $\al_n=1-\al n ^{-1}$ with  $\al >2$ (for the momentum coefficient), along with $\ld_n =\ld  \in (0,2\beta)$ (for the proximal indexes), and  $w_n=1-\rho n^{-2}$  and
with  $0< \rho < \al (\al -2) ( 1- \frac{\ld}{4 \beta} )$ (for the relaxation factors).
In this  framework (see \cite[corallary 4.9]{ac}, they have reached the   estimates  $\|\xnp-\xn\|^2= {\cal O}(n ^{-1})$ and  $\sum_{n } n\|\xnp-\xn\|^2 < \li$ (for the discrete velocity),
 along with     $\sum_{n} n^{-1}  \|G^I_{\ld}(x_n)\|^2< \li$ and  $\lim_{n \to \li} \|G^I_{\ld}(x_n)\|=0$  (for the  fixed point residuals). \\
  \\
%
\ind Thanks to the correction term in CRIFBA we improve these  last rates  from $\|\xnp-\xn\|^2= {\cal O}(n ^{-1})$ and
 $\|G^I_{\ld}(x_n)\|^2={\cal O}(1)$  to $\|\xnp-\xn\|^2= o(n ^{-2})$ and
 $\|G^I_{\ld}(x_n)\|^2=  o(n ^{-2})$.
%
 \subsection{Organization of the paper.} An outline of this paper is as follows. In section \ref{pado2}, we present CRIFBA (with full details on the parameters)   and  its main convergence results. A proof of the main results is then proposed.  In section \ref{pado3}, we specialize  CRIFBA to the setting of the more general monotone inclusion (\ref{gmoni}). An application of  CRIFBA is  given in section \ref{pado4} relative to convex-concave saddle point problems. \\
%
  \begin{remark}{From now on,  so as to simplify the notations,   we   (often) use the  following notation: given any sequence $(u_n)$, we denote  $\dot{u}_n=u_n-u_{n-1}$. \\
 } \end{remark}
\section{\label{pado2} \large Main results and preliminary estimations.}
\setcounter{equation}{0}
 \subsection{Main convergence results.}
  The following result states  the convergence of CRIFBA,   with an     accuracy    measured   through  the operator $G_{\ld}^M$ (introduced in (\ref{defmld})). \\
 \begin{theorem} \label{thmg}
 Let $L, M: \Hs \to \Hs$  be linear self-adjoint and positive definite maps, and let $A: \Hs \to 2^{\Hs}$ and $B :\Hs \to \Hs$ verify (\ref{cdi}),  with     $S:=(A+B)^{-1}(0) \neq \emptyset$.
Suppose that   $\{x_n, v_n \} \subset \Hs $ are    generated by CRIFBA  with   $\{s_1, s_0, \nu_0\} \subset[0,\li)$  and
$\{e,   \ld, w\}\subset(0,\li)$  verifying
 \begin{eqnarray} \label{gfsc}
   \label{gfsc1}
   && \mbox{$ 2 s_1 <    s_0 <   e$}, \es
  \mbox{$0 < w <  1$}.
    \end{eqnarray}
 Suppose furthermore     that one of the following conditions (\ref{cpal1}) and  (\ref{cpal2}) holds :
  \begin{subeqnarray}
   \label{cpal}
  \label{cpal1}
 &&    \mbox{$\exists \delta>0$ s.t. $ \ld  \|L \| \le  4 \delta$  \es  and \es  $\bar{M}_1:=M- \frac{\delta}{w(1- w)} I$ is positive definite  }, \\
 \label{cpal2}
 &&  \mbox{$\bar{M}_2:= M- \frac{\ld }{w(1- w)} L$ is positive definite}.
 \end{subeqnarray}
Then the  following  properties  are reached:
   \begin{subeqnarray}
   \label{borg}
    \label{borg1}
   &  &     \mbox{ $\| \dxnp  \|_M^2 =o( n^{-2})$,
   \es     $\sum_n n  \| \dxnp  \|_M^2 <  \li$, \es $\sum_n n^2  \| \dxnp -\dxn  \|_M^2 <  \li$}, \\
    \label{borg2}
   &  &    \mbox{ $ \|v_n \|_M^2 =o( n^{-2})$,
   \es $\sum_n n   \| v_n\|_M^2 <  \li$, \es $\sum_n n ^2  \| \dvnp \|_M^2 <  \li$}, \\
    \label{borg3a}
    &  &    \mbox{ $ \|G_{\ld}^M (x_n) \|_M^2 =o( n^{-2})$,
   \es $\sum_n n   \|G_{\ld}^M (x_n) \|_M^2 <  \li$}. \\
      &  & \hskip -2.cm \mbox{Moreover,  denoting  $\yn= x_n + \lb 1-\frac{1}{ w}\rb v_{n}$, we have:} \nonumber \\
 \label{borg1b}
   &  &     \mbox{ $\| \dynp  \|_M^2 =o( n^{-2})$,
   \es     $\sum_n n  \| \dynp  \|_M^2 <  \li$}, \\
     \label{borg3c}
    &  &    \mbox{ $ \|G_{\ld}^M (y_n) \|_M^2 =o( n^{-2})$,
   \es $\sum_n n   \|G_{\ld}^M (y_n) \|_M^2 <  \li$}. \\    
      &  & \hskip -2.cm \mbox{If, in addition, $M$ and $L$ are bounded, then:} \nonumber \\
 \label{borg5}
   &  &  \mbox{$\exists \bar{x}\in S$, s.t. $(x_n,y_n)  \rightharpoondown (\bar{x},\bar{x}) $ weakly in $(\Hs, \|.\|_M)^2$}, \\
      \label{borg6}
   &  &  \mbox{$\exists y_n^* \in (A+B)(y_n)$, s.t. $\|y_n^*\|_M=o( n^{-2})$ and $\sum_n n\|y_n^*\|_M^2< \li$}.
   \end{subeqnarray}
\end{theorem}
\ind  Theorem \ref{thmg}  will be proved in Section \ref{pthmg}. \\
\\
\ind Let us give some  comments on the above theorem.
\begin{remark} \label{remg} { Recall that a linear self-adjoint operator ${\cal M}: \Hs \to \Hs$ is called positive definite if it satisfies  $\inf_{x \in \Hs} \frac{\|\la{\cal M}x, x\ra\|}{ \|x\|} >0$.  So, we  emphasize,  for a bounded operator  $L$,  that the two  conditions    (\ref{cpal1}) and   (\ref{cpal2})  are always feasible for $\delta$ and $\ld$ small enough. In specific,  regarding the classical setting of (\ref{pbg})-(\ref{cdi}) when $L=\beta^{-1} I$ (with $\beta >0$) and $M=I$ (as used  for instance in
\cite{ac,rafape,moudol}),  (\ref{cpal1}) reduces to $0< \ld  \le  4 \beta \delta $ and  $0 <  \delta < w(1-w)$, for some $\delta >0$,
which is    equivalent to  $0< \ld  <  4 \beta w(1-w)$,
    while    (\ref{cpal2}) becomes the more stringent condition $0 < \ld  < w(1-w)\beta$. Nonetherless,  condition  (\ref{cpal2}) (in general) does not require $L$ to be bounded.
 }\end{remark}
%
 %
 \subsection{Preliminary estimations on  CRIFBA.} In order to prove Theorem \ref{thmg}, we  exhibit  a Lyapunov sequence in connection with the proposed algorithm. This  allows us to obtain a series of preliminary estimates. Next, we derive additional estimates from a suitable reformulation  of CRIFBA in terms of the quantities $\dxn$ and $v_n$. Finally we combine the two series of results so as to reach the desired estimates.  A preliminary  observation regarding this section is given by  the following remark. 
 \begin{remark} \label{dax} {It is importance to  notice that (\ref{sc1b}) can be reformulated as
 \be \label{nsc1b} \xnp = z_n-w \ld G_{\ld}^M(z_n), \ee
 where $G_{\ld}^M$ was introduced in (\ref{defmld}). \\
 }\end{remark}
  \subsubsection{A useful  reformulation of the algorithm.} 
  We begin with providing  a useful  reformulation of CRIFBA.  As standing assumptions we assume that  $L: \Hs \to \Hs$ and $M: \Hs \to \Hs$  are linear self-adjoint and positive definite maps, and that  $A: \Hs \to 2^{\Hs}$ and $B :\Hs \to \Hs$ verify condition (\ref{cdi}).  A key result in our analysis is  given by   the following proposition.\\
\begin{proposition}\label{jeu0}   The iterates
     $(x_n)$ and $(v_n)$ generated by CRIFBA satisfy  (for $n \ge 0$)
     \begin{subeqnarray}
\label{subt}
\label{subt1}
 & & \mbox{$v_{n+1}= (\ld w)  G_{\ld}^M(z_n)$,}\\
 \label{subt2}
 & &   \dxnp    + v_{n+1}   = \ttn \dxn + \gamma_{n}  v_{n}.
   \end{subeqnarray}
  \end{proposition}
\ddem  For $n \ge 0$, we have  $v_{n+1}=z_{n}-  \xnp  $ (by definition of $v_n$)  together with  $z_n=\xnp + (\ld w)  G_{\ld}^M (z_n) $ (from (\ref{nsc1b})), hence,    we obviously infer that
  $v_{n+1} = (\ld w) G_{\ld}^M(z_{n})$, which yields (\ref{subt1}). Furthermore, by    $z_n= x_n + \ttn   \dxn  + \gamma_n  v_n $ (from (\ref{sc1c})), we immediately obtain 
   $v_{n+1}:=z_n-\xnp=\ttn   \dxn  +\gamma_n  v_n  -  \dxnp$, which   entails (\ref{subt2})\edem 
\subsubsection{Co-coerciveness of $G_{\ld}^M$ and basic properties on $(x_n)$ and  $(v_n)$.}
 The following result  establishes  a  co-coerciveness property for   $G_{\ld}^M$, which   plays a central  role in our methodology. \\
   \begin{proposition} \label{gcd} For any  $\ld >0$   and for any $(x_1,x_2) \in \Hs^2$,  we have
    \be \label{mohg} \ba \hskip 0.8cm  \la \Delta G _{\ld}^M(x_1,x_2)  , x_1-x_2 \ra _M \\
   \hskip 2.cm  \ge     \| \Delta B (x_1,x_2) \|^2_{L^{-1}}   +  \ld  \|   {B}(x_1)-{B}(x_2) \|_M^2
    -  \ld \la   {B}(x_1)-{B}(x_2) ,\Delta B (x_1,x_2) \ra, \ea
 \ee
 where $\Delta G _{\ld}^M(x_1,x_2)=G_{\ld}^M(x_1)- G_{\ld}^M(x_2)$ and   $\Delta B (x_1,x_2)=B(x_1)- B(x_2)$. So, for $i=1,2$, the following inequalities hold:
\be
 \label{mohga}
 \hskip 1.cm \la \Delta {G}_{\ld}^M(x_1,x_2) , x_1-x_2 \ra _M \ge   \al_i     \| \Delta {B}(x_1,x_2)\|^2_{L^{-1}}
     +   \ld \la  H_i \Delta {G}_{\ld}^M(x_1,x_2), \Delta {G}_{\ld}^M(x_1,x_2) \ra ,
 \ee
 with parameters  $\al_i$ and $H_i$ defined by
 \begin{subeqnarray}
 \label{dfal}
 \label{dfal1}
&& \mbox{$\al_1=   1 -  \frac{\ld}{4 \delta }  \|L \|$ \es and \es $H_1=  M -\delta I$ (for some given $\delta>0$)}, \\
\label{dfa2}
&& \mbox{$\al_2= \frac{3}{4  }$ \es and \es $H_2=  M -\ld  L$}.
 \end{subeqnarray}
  \end{proposition}
\ind The proof of Proposition \ref{gcd} is given in Appendix \ref{pgcd} and it makes use of  the following observation. \\
\begin{remark}\label{pop} {Recall that for any linear self-adjoint and positive definite operator ${\cal M}: \Hs \to \Hs$, there exists a linear self-adjoint and positive definite,  map denoted ${\cal M}^{\frac{1}{2}} $ such that  ${\cal M}={\cal M}^{\frac{1}{2}}{\cal M}^{\frac{1}{2}}$ (see \cite{sebtarc}). The operator ${\cal M}^{\frac{1}{2}}$ is called of roof of ${\cal M}$ and it also satisfies
  $\|{\cal M}^{\frac{1}{2}} \|^2=\|{\cal M} \|$. \\
  } \end{remark}
\\
\ind As an immediate consequence of the previous proposition,  we provide basic properties concerning  the sequences $(x_n)$ and $(v_n)$. \\
 \begin{proposition} \label{rilo}
   Suppose that   $\{x_n, v_n\} \subset \Hs$  are   generated by   CRIFBA  with  parameters $\{  s_0, s_1 \} \subset [0,\li)$,   $\{  e,\nu_0,  \ld  \}\subset (0,\li)$ and  $w \in (0,1)$.
 Suppose furthermore, for  some integer  $i_c\in \{1,2\}$,     that $\bar{M}_{i_c}$ is positive definite,
 where $\bar{M}_{1}$ and  $\bar{M}_{2}$ are defined by
  \be
   \label{cout}
  \mbox{$\bar{M}_1:= M- \frac{\delta}{w(1- w)} I$ is   (for some given $\delta>0$)}, \es \es
  \mbox{$\bar{M}_2:= M- \frac{\ld }{w(1- w)} L$}.
 \ee
 Then,  for  $n \ge 1$ and for  any $q \in S$ we have the following inequalities:
 \begin{subeqnarray}
   \label{outg}
  \label{out1}
  && \hskip 1.cm \mbox{$\la v_n, x_n-q \ra_M \ge (\ld w) \al_{i_c}      \|  {B}(z_{n-1}) - {B}(q) \|^2_{L^{-1}} $} +   \mbox{ $\frac{(1- w)^2}{w}  \| v_n\|_{ M } ^2$}, \\
 \label{out2}
 && \hskip 1.cm  \mbox{$\la \dvnp , \dxnp \ra_M \ge (\ld w)  \al_{i_c}    \|  {B}(z_{n}) -  {B}(z_{n-1}) \|^2_{L^{-1}} +  \frac{(1- w)^2}{w} \| \dvnp\|_{ M } ^2$},
 \end{subeqnarray}
 where $\al_{i_c}$ is given by (\ref{dfal}). \\ \\
 \end{proposition}
 \ddem For convenience of the reader  we set $\zeta_n=G_{\ld}^M(z_{n})$.  Clearly, for $n \ge 1$,  by Remark \ref{dax} we have
  $x_n= z_{n-1}- \ld w \zeta_{n-1}$.
  It follows immediately for $q \in S$  that  \\
  \ind $ \ba  \la v_n , x_n-q \ra_M =  (\ld w) \la \zeta_{n-1}, x_n-q \ra_M \\
  \hskip 2.cm = (\ld w) \la \zeta_{n-1}, z_{n-1}-q \ra_M - (\ld w)^2 \| \zeta_{n-1}\|^2_M. \ea$ \\
  Moreover,  by Proposition \ref{gcd} and recalling that $G_{\ld}^M(q)=0$ we have  (for $i=1,2$)\\
  \ind $ \la \zeta_{n-1}, z_{n-1}-q \ra_M \ge \al_i    \|  {B}(z_{n-1}) - {B}(q) \|^2_{L^{-1}}
  + \ld  \la H_i  \zeta_{n-1},  \zeta_{n-1} \ra $.\\
 Consequently, by the previous  statements,    we  deduce that \\
  \ind $\ba \la v_n , x_n-q \ra_M - (\ld w) \al_i    \|  {B}(z_{n-1}) - {B}(q) \|^2_{L^{-1}} \\
   \hskip 2.cm \ge    \la \lb  \ld^2 w H_i   -(\ld w)^2 M ) \rb  \zeta_{n-1}   , \zeta_{n-1}\ra \\
    \hskip 2.cm  =  \ld^2 w\la \lb  H_i- w M   \rb  \zeta_{n-1}   , \zeta_{n-1}) \ra =  \ld^2 w  (1- w)\la   \bar{J}_i  \zeta_{n-1}, \zeta_{n-1}\ra,   \ea $ \\
  where \\
  \ind $\ba \bar{J}_1 = (1- w)^{-1}(H_1- w M)=  (1- w)^{-1}(M-\delta I - w M) \\
  \ind = M- \delta (1- w)^{-1}I= (1-w) M + w \bar{M}_1, \ea $ \\
  and \\
   \ind $\ba \bar{J}_2 = (1- w)^{-1}(H_2- w M)=  (1- w)^{-1}(M-\ld  L - w M)\\
   \ind = M- \ld (1- w)^{-1}L = (1-w) M + w \bar{M}_2.  \ea $ \\
 This  leads  immediately  to (\ref{out1}).  Next,  for $n \ge 1$, by (\ref{subt1}) we obviously have
 $\ba  \dvnp = \ld w  (  \zeta_{n} -  \zeta_{n-1}  ), \ea $
 while Remark \ref{dax} yields  $ \dxnp= \dzn -  \ld w (  \zeta_{n} -  \zeta_{n-1}  )$.
 Then  we   immediately see   that \\
  \ind $ \la \dvnp, \dxnp \ra_M
 =  \ld w \la  \zeta_{n} -  \zeta_{n-1} , \dot{z}_{n} \ra_M -(\ld w)^2\|  \zeta_{n} -  \zeta_{n-1}  \|^2_M.$ \\
 In addition, by  Proposition \ref{gcd}, we have  (for $i=1,2$)\\
 \ind $ \ba  \la \dot{\zeta}_{n} , \dzn \ra_M \ge  \al_i    \|  {B}(z_{n}) -  {B}(z_{n-1}) \|^2_{L^{-1}}  +  \ld \la H_i  \dot{\zeta}_{n} ,  \dot{\zeta}_{n}  \ra.  \ea $ \\
  Then by the previous two results we are led to \\
   \ind $\ba \la \dvnp, \dxnp \ra_M - (\ld w) \al_i    \|  {B}(z_{n}) -  {B}(z_{n-1}) \|^2_{L^{-1}} \\
  \hskip 2.cm  \ge
     \la ( \ld^2w H_i - (\ld w)^2 M)  \dot{\zeta}_{n} , \dot{\zeta}_{n}  \ra  =  \ld^2w \la (  H_i - w M)  \dot{\zeta}_{n} , \dot{\zeta}_{n} \ra \\
 \hskip 2.cm =     \ld^2w (1-w) \la  \bar{J}_i  \dot{\zeta}_{n} , \dot{\zeta}_{n}  \ra
 ,   \ea $ \\
where  $\bar{J}_i$ was   previously introduced.   This  readily amounts to (\ref{out2})\edem 
 \subsubsection{Links between  the iterates and the graph of  $(A+B)$.} 
  Consider  the  elements $y_n$ and $y_n^*$ defined by
 \begin{eqnarray}
 \label{diss}
   \mbox{$\yn= -\frac{1}{ w}v_{n} + z_{n-1}  $  and $\yn^*= (\ld w)^{-1} Mv_n + B(\yn)-  B(z_{n-1})  $}.
 \end{eqnarray}
 \begin{remark}{ Note that, for $n\ge 1$ we have  $z_{n-1}=v_n+x_n$ (from definition of $v_n$), so $y_n= x_n+ (1-\frac{1}{ w}) v_n$, which is nothing but the formulation of $y_n$ used in Theorem \ref{thmg}.\\ \\
  } \end{remark}
  \ind The following result  makes  the  connection between  $\yn^*$ and $(A+B) (\yn)$. \\
 \begin{proposition} \label{jeu1} Let
 $\{ \xnp, v_n \}_{n\ge 0}$ be sequences   produced by CRIFBA. Then, for $n\ge 0$, the elements $\ynp$ and  $v_{n+1}$ given by
 (\ref{diss}) satisfy
 \begin{subeqnarray}
 \label{piss}
 \label{piss0}
 &  &  \mbox{$  (\ld  w)^{-1} M v_{n+1} \in   B(z_n) +  A (\ynp) $}, \\
 \label{piss1}
 &  &  \mbox{$ \ynp^*   \in   B (\ynp) +  A (\ynp) $}.
 \end{subeqnarray}
 If, in addition, $M$ and $L$ are bounded and  such that  $M-\rho L$ or  $M-\rho I$ is positive definite (for some $\rho >0$),  then
\be  \label{piss2}
   \mbox{$\|\ynp^*\|_M \le \frac{1}{ w}\lb   \frac{1}{ \ld } \|M\| + \frac{1}{ \rho^{1/2} } (\|M\|.\|L\|) ^{\frac{1}{2}}
    (1+\|L\|^{\frac{1}{2}})   \rb  \|\vnp\|_M $}.
    \ee
 \end{proposition}
 \ind The proof of Proposition \ref{jeu1} is given in Appendix \ref{djeu1}. \\
 %
\section {\large Convergence analysis of CRIFBA.}
\setcounter{equation}{0}
 %
 A series of estimates are obtained here  by means of a Lyapunov analysis (based upon Proposition \ref{gfru0}) and using the reformulation of CRIFBA. The main results of Theorem \ref{thmg} will be derived  as a combination of the previous series of  estimates. As standing assumptions we assume that  $L: \Hs \to \Hs$ and $M: \Hs \to \Hs$ are  linear self-adjoint and positive definite maps and that $A: \Hs \to 2^{\Hs}$ and $B :\Hs \to \Hs$ are  maximally monotone operators.
 \subsection{Estimates from an energy-like sequence.}
   With the iterates $\{x_n,v_n\}$ produced by CRIFBA, we associate the sequence  $( {\cal E}_n(s,q))$ defined   for $(s,q) \in (0, \li) \mul  S$ and for $n \ge 0$ by
  \be \label{lf}  \ba {\cal E}_n(s,q)=  \frac{1}{2}    \| s (q-\xn)- \nu_n \dxn  \|^2_M +  \frac{1}{2}  s(e  -s) \|\xn -q\|^2_M 
  +  s   (e+\nu_{n})  \la v_n,\xn-q \ra_M. \ea \ee
\ind Our Luapunov analysis we be based upon the following lemma. \\
  \begin{lemma} \label{estimg0}
 Suppose that (\ref{cdi}) holds   and   that  $\{x_n, v_n\} \subset \Hs$  are     generated by   CRIFBA  with  parameters $\{  s_0, s_1 \} \subset [0,\li)$ and  $\{  e,\nu_0, w, \ld  \}\subset (0,\li)$  such that
   \be \label{fsc}
    \mbox{$0 \le  s_1 <   s_0 <e$}, \es \es
   \mbox{$0 < w < 1$}.
   \ee
Suppose furthermore that there exists   $i_c\in \{1,2\}$ such that the following condition holds:
\be \label{cic} \mbox{ $\al_{i_c}$ is non-negative  and  $\bar{M}_{i_c}$ is positive definite}, \ee where  $\al_{i_c}$ and  $\bar{M}_{i_c}$ are given by (\ref{dfal}) and (\ref{cout}), respectively.
Then, for any $(s,q)\in (0,e] \mul \Hs$ and for $n \ge 1$, we have
\be  \label{dju2}
 \ba \dot{\cal E}_{n+1}(s,q) + s (s_0-s_1) \la v_n, \xn -q \ra_M \\
 + \frac{(1-w)}{w}  (e+\nu_{n+1}) (e-s +\nu_{n+1}) \| \dvnp \|^2_{{M}}  \\
 +  \frac{1}{2} (e+\nu_{n+1})\| \dxnp-\ttn \dxn \|^2 _M \\
     +   \lb s_0  -s \rb (e+\nu_{n+1}) \la \vn , \dxnp \ra_M
       +    \frac{1}{2}\lb e  -s \rb ( e+2  \nu_{n+1} )  \|\dxnp\|^2 _M\le 0. \ea \ee
  In particular, for $q \in S$,   the sequence $({\cal E}_n(s_0,q))_{n \ge 1}$ is non-increasing  and  convergent,  and  the  following estimates are reached:
\begin{subeqnarray}
\label{pp}
\label{pp0}
&& \mbox{$\sup_{n \ge 1}  \|x_n-q\|^2 _M\le \frac{2{\cal E}_{1}(s_0,q) }{s_0(e- s_0)}$}, \\
\label{pp1}
&& \mbox{$\sup_ {n \ge 1} \nu_n  \la v_n  , \xn-q \ra_M \le \frac{{\cal E}_{1}(s_0,q)}{ s_0  }$}, \\
\label{pp2}
&&   \mbox{$\sup_n n \| \dxn \|_M < \li$}, \\
\label{pp3}
&& \mbox{$ \sum_{n \ge 1} \nu_{n+1}^2 \|\dvnp \|^2_{ {M} } \le  \frac{w {\cal E}_{1}(s_0,q)}{(1-w)^2}$}, \\
\label{pp4}
&& \mbox{$\sum_{n \ge 1}       \nu_{n+1}   \|\dxnp\|^2_M \le \frac{ {\cal E}_{1}(s_0,q) }{e-s_0} $}, \\
\label{pp5}
&&   \mbox{$\sum_{n \ge 1}  \la v_n, \xn-q \ra_M \le  \frac{{\cal E}_{1}(s_0,q)}{s_0 (s_0-s_1)} $}, \\
\label{pp6}
&& \mbox{$\sum_{n }     n^2  \|\dxnp-\dxn\|^2 _M< \li$}.
\end{subeqnarray}
  \end{lemma}
 \subsubsection {Proof of Lemma \ref{estimg0}.}
  \paragraph{Proof of Lemma \ref{estimg0} - Part I : a useful equality  for a Lyapunov analysis.}
An  important equality of independent interest is  proposed   here  relative  to our  method through the wider  framework of   sequences   $\{x_n, d_{n}\} \subset   \Hs$ and   $\{e, \nu_n, \ttn \} \subset  (0,\li)$  verifying (for $n \ge 0$)
\begin{subeqnarray}
\label{adr0}
\label{adr0a}
&  &   \dxnp   - \ttn \dxn +  d_n   =0, \\
\label{fdc}
&  &   (e+ \nu_ {n+1}) \ttn= \nu_n.
\end{subeqnarray}

\ind  As a key element of our methodology  we associate with (\ref{adr0})  the quantity  $F_n(s,q)$ given  for any  $(s,q)\in [0,\li) \mul \Hs$ by
\be
\label{ener}
 \mbox{$F_n(s,q) =  \frac{1}{2}   \| s (q-\xn)-  \nu_n \dxn    \|^2_M+   \frac{1}{2} s(e  -s) \|\xn -q\|^2_M$.}
 \ee
 \ind Basic properties  regarding
 the sequence  $({F}_{n}(s,q))$  are established through  the following proposition.\\
\begin{proposition} \label{gfru0}  Let    $\{x_n, d_n \} \subset \Hs$ and $\{  \ttn, \nu_n, e \} \subset (0,\li)$ verify (\ref{adr0}),
 and suppose that  $M: \Hs \to \Hs$  is  a linear self-adjoint and positive definite map.
 Then for   $(s,q)  \in  \lb 0, e \right ] \mul  \Hs  $   and for $n \ge 0$ we have
 \be \label{grfd0}
   \ba \dot{F}_{n+1}(s,q) + \frac{1}{2} (e+\nu_{n+1}) ^2 \| \dxnp- \ttn \dxn \|^2_M \\
   +   s  (e+\nu_{n+1}) \la d_{n}, \xnp -q \ra_M   \\
   + (e+\nu_{n+1}) (e-s+\nu_{n+1})\la  d_{n} ,\dxnp  \ra_M
 = -  \frac{1}{2}\lb e  -s \rb ( e+2  \nu_{n+1} )  \|\dxnp\|^2_M  .
    \ea   \ee
 \end{proposition}
\ind The proof of Proposition \ref{gfru0} is given in Appendix \ref{pap2}. \\
   \paragraph{Proof of Lemma \ref{estimg0} - Part II.}
  Let us begin with proving (\ref{dju2}). It can be observed (from   (\ref{subt2})) that   the iterates  $\{x_n, v_n\}$ generated by CRIFBA enter the special case of  the general iterative process  (\ref{adr0}) when taking
  $d_n= v_{n+1}-\gamma_n v_n $. Hence, for $n \ge 0$,  by Proposition \ref{gfru0} we obtain
  \be \label{ffd0}
  \hskip .5cm  \ba \dot{F}_{n+1}(s,q) + \frac{1}{2} \tau_n^2 \|W_n\|^2_M \\
  +  (s \tau_n) \la d_n, \xnp -q \ra _M + (\tau_n^2 \gamos)\la  d_n ,\dxnp  \ra_M       =   -  \frac{1}{2}\lb e  -s \rb ( e+2  \nu_{n+1} )  \|\dxnp\|^2_M .
    \ea   \ee
    where $\tau_n=e+ \nu_{n+1}$,   $\gamos=1- s \tau_n^{-1}$ and   $W_n=\dxnp-\ttn \dxn$. Let us  evaluate  the quantity  $(s \tau_n) \la d_n, \xnp -q \ra _M + (\tau_n^2 \gamos)\la  d_n ,\dxnp  \ra_M$. 
      Setting $U_n=  \la v_{n}, \xn-q \ra_M$,  by    $d_n= v_{n+1}-\gamma_n v_n $ we have
    \begin{eqnarray}
    && \mbox{$\la d_n, \xnp -q \ra_M = U_{n+1}  + \la -\gamma_n v_n, \dxnp \ra_M   -\gamma_n U_n $} ,\\
   && \mbox{$\la d_n, \dxnp \ra_M = \la v_{n+1},\dxnp \ra_M  + \la -\gamma_n v_n, \dxnp \ra_M$}.
   \end{eqnarray}
     This, noticing  that $ s\tau_n + \tau_n^2 \gamos=\tau_n^2$, amounts to
     \[\ba  (s \tau_n) \la d_n, \xnp -q \ra _M + (\tau_n^2 \gamos)\la  d_n ,\dxnp  \ra_M \\
     \hskip 2.cm = (s\tau_n)  U_{n+1} -  (s\tau_n)\gamma_n U_n + (\tau_n^2 \gamos)  \la v_{n+1},\dxnp \ra_M +  \tau_n^2\la -\gamma_n v_n, \dxnp \ra _M . \ea   \]
    In addition, we obviously have \\
 \ind  $  \la v_{n+1},\dxnp \ra_M=  \la \dot{v}_{n+1},\dxnp \ra_M+ \la v_n,\dxnp \ra_M$. \\
 Then, by the previous arguments  we obtain  \\
 \ind $\ba   (s \tau_n) \la d_n, \xnp -q \ra _M + (\tau_n^2 \gamos)\la  d_n ,\dxnp  \ra_M \\
 \ind \ind= (s\tau_n)  U_{n+1}-\gamma_{n}  (s\tau_n)U_n +(\tau_n^2 \gamos) \la  v_{n+1},\dxnp \ra_M -\gamma_{n} \tau_n^2 \la v_n, \dxnp \ra_M \\
 \ind \ind = (s\tau_n)  U_{n+1}-\gamma_{n}  (s\tau_n)U_n  +(\tau_n^2 \gamos) \la  \dvnp,\dxnp \ra_M+ \tau_n^2( \gamos-\gamma_{n} ) \la v_n, \dxnp \ra_M
    . \\
    \ea$ \\
   Furthermore, observing   that $\gamma_n=1-s_0 \tau_n^{-1}$ and that $\dot{\tau}_n :={\tau}_n-{\tau}_{n-1}  =s_1$ , an easy computation  gives us  \\
   \ind $\gamma_{n}  \tau_n =\tau_{n}-s_0 =  \tau_{n-1} - (s_0-s_1),  $\\
   \ind $ \tau_n ( \gamos-\gamma_{n} )= (s_0-s).  $\\
 In addition, for $n \ge 1$, by Proposition \ref{rilo} and setting  $\eta=\frac{(1-w)^2}{w}$, we have \\
  \ind $\la \dvnp, \xnp \ra_M \ge  \eta  \|\dvnp\|^2 _{{M}}$. \\
 Therefore combining the previous results  amounts to
 \be \ba    (s \tau_n) \la d_n, \xnp -q \ra _M + (\tau_n^2 \gamos)\la  d_n ,\dxnp  \ra_M \\
   \hskip 2.5cm\ge
       (s\tau_n)  U_{n+1} -  (s\tau_{n-1})  U_{n}  + s (s_0-s_1)U_n \\
    \hskip 2.5cm    + \eta (\tau_n^2  \gamos)  \|\dvnp\|^2_{{M}}  +  \tau_n  (s_0-s)  \la v_n, \dxnp \ra_M
        . \ea \ee
 Consequently,   in light of (\ref{ffd0}), we deduce for $n \ge 1$ that
   \[
    \ba \dot{F}_{n+1}(s,q)  +   (s \tau_n) U_{n+1}-    (s \tau_{n-1}) U_{n}  + s (s_0-s_1)U_n \\
   \hskip 0.5cm
     +  \eta \tau_n^2  \gamos  \|\dvnp\|^2 _{{M}} +  \tau_n  \lb s_0 -s \rb  \la \vn , \dxnp \ra_M
       \le   -   \frac{1}{2}\lb e  -s \rb ( e+2  \nu_{n+1} )  \|\dxnp\|^2_M .
    \ea   \]
 This, from  ${\cal E}_{n}(s,q)={ F}_{n}(s,q) +   (s \tau_{n-1}) U_{n} $  (in light of (\ref{lf})),  can be rewritten as
   (\ref{dju2}). \\
    \\
    \ind Now we prove the second part of Lemma \ref{estimg0}.    For $q\in S$ and $s=s_0$,  inequality (\ref{dju2}) becomes
    (for $n \ge 1$)
  \be \label{dju2b}
  \ba \dot{\cal E}_{n+1}(s_0,q) + s_0 (s_0-s_1) \la v_n, \xn-q \ra _M \\
  \hskip 2.cm
     +  \eta \tau_n ^2 \gamma_n     \| \dvnp \|^2_{{M}}  +  \frac{1}{2} \tau_n^2\|W_n\|^2_M  +   \frac{1}{2}\lb e  -s_0 \rb ( e+2  \nu_{n+1} )  \|\dxnp\|^2_M \le 0.
    \ea   \ee
   Clearly, we know that $U_n$ is non-negative (from Proposition \ref{rilo}). It follows immediately that the non-negative sequence  $({\cal E}_{n}(s_0,q))_{n \ge 1}$ is non-increasing, since the constants  $s_0-s_1$ and $e-s_0$ are non-negative (in light of condition (\ref{fsc})). Whence  it is convergent and bounded. Also recall  from  (\ref{lf}) that
    \be \label{lf2} \es \es  \ba {\cal E}_n(s_0,q)=  \frac{1}{2}   \| s_0  (q-\xn)- \nu_n \dxn  \|^2_M  +  \frac{1}{2}  s_0 (e-s_0) \|\xn -q\|^2 _M    +  s_0( e+2  \nu_{n} ) \la v_n, \xn-q \ra _M
  . \ea \ee
   Then, for $n \ge 1$, by  the inequality ${\cal E}_{n}(s_0,q)\le {\cal E}_{1}(s_0,q)$  we get
   \begin{eqnarray}
   &&   \mbox{$ \frac{1}{2} s_0 (e-s_0) \|\xn -q\|^2 _M \le  {\cal E}_{1}(s_0,q)$},\\
   &&  \mbox{$ s_0  ( e+2  \nu_{n} ) \la v_n, \xn-q \ra _M  \le  {\cal E}_{1}(s_0,q)$}, \\
   &&  \mbox{$\| \nu_n \dxn  \|_M- \| s_0  (q-\xn)\|_M  \le \sqrt{ {\cal E}_{1}(s_0,q)}$}.
    \end{eqnarray}
    Estimates (\ref{pp0}), (\ref{pp1}) and (\ref{pp2})  are direct consequences of these last three inequalities.
     Moreover, given some integer $N \ge 1$, by adding  (\ref{dju2b}) from $n=1$ to $n=N$  we obtain
   \be
  \ba {\cal E}_{N+1}(s_0,q)     +  s_0 (s_0-s_1)\sum_{n=1}^N \la v_n, \xn-q \ra _M  \\
   +     \eta   \sum_{n=1}^N  \tau_n ^2 \gamma_n \| \dvnp \|^2_{{M}}  +  \frac{1}{2} \sum_{n=1}^N \tau_n^2 \|\dxnp-\ttn \dxn\|^2_M \\
   + \frac{1}{2} \lb e  -s_0 \rb  \sum_{n=1}^N  ( e+2  \nu_{n+1} )  \|\dxnp\|^2_M       \le  {\cal E}_{1}(s_0,q).
    \ea   \ee
    It follows immediately that  
    \begin{subeqnarray}
    \label{joi}
     \label{joi1}
   &&  \mbox{$  \eta \sum_{n=1}^N  \tau_n ^2 \gamma_n \| \dvnp \|^2_{{M}}  \le  {\cal E}_{1}(s_0,q) $}, \\
    \label{joi2}
   &&  \mbox{$\frac{1}{2} \lb e  -s_0 \rb  \sum_{n=1}^N  ( e+2  \nu_{n+1} )  \|\dxnp\|^2_M   \le  {\cal E}_{1}(s_0,q)$}, \\
    \label{joi3}
   &&   \mbox{$s_0 (s_0-s_1)\sum_{n=1}^N \la v_n, \xn-q \ra _M   \le  {\cal E}_{1}(s_0,q)$}, \\
    \label{joi4}
   &&   \mbox{$\frac{1}{2} \sum_{n=1}^N \tau_n^2 \|\dxnp-\ttn \dxn\|^2_M  \le {\cal E}_{1}(s_0,q)$}.
    \end{subeqnarray}
    This straightforwardly yields (\ref{pp3}), (\ref{pp4}) and (\ref{pp5}). The last estimate (\ref{pp6}) is simply deduced  from
  (\ref{joi4}) and  (\ref{pp4}) (in light of the definition of $\ttn$)\edem \\
    \subsection{Estimates from the reformulation of the method.} Additional estimates are established regarding CRIFBA, especially on the sequence $(v_n+ \dxn)$. \\
  \begin{lemma} \label{estimg2} Let     $\{x_n, v_n\} \subset \Hs$  be  given by   CRIFBA with  $\{s_1, s_0\} \subset [0,\li)$ and $\{ e, \nu_0,  w, \ld  \}\subset (0,\li)$  verifying $0 < s_0 <  e$.
Then, for any $(s,q)\in (0,e] \mul  \Hs$ and for $n \ge 1$, we have
  \be  \label{mst} \ba (e+\nu_{n+1})^2\| v_{n+1}+ \dxnp \|^2 _M- (e+\nu_{n})^2 \|  v_{n}+ \dxn \|^2_M  \\
  \hskip 1.cm +   (s_0-  2 s_1) (e+\nu_{n+1})  \|  v_{n}+ \dxn \|^2_M    \le   s_0^{-1} (e-s_0+s_1)^2 (e+\nu_{n+1})  \|\dxn\|^2_M . \ea \ee
 Suppose, in addition,  that condition (\ref{cic}) holds and that the parameters verify
  \be  \label{mstc} \mbox{$0 \le s_1 < (1/2) s_0$ \es and \es $0<w<1$}. \ee
   Then the  following estimates are reached
   \begin{subeqnarray}
\label{dd}
\label{dd1}
&& \mbox{$\sum_n n   \|\dxn+ v_n \|^2_M < \li$}, \\
\label{dd2}
&& \mbox{$\|\dxn+ v_n\|_M=o(n^{-1})$}.
\end{subeqnarray}
   \end{lemma}
\ddem
  For $n\ge 1$, according to Lemma \ref{jeu0} and denoting $\tau_n =e+\nu_{n+1}$ we have  $v_{n+1}+ \dxnp =  \ttn \dxn + \gamma_n v_n$,
 together with $\gamma_n= 1-s_0 \tau_n^{-1}$ and $\ttn = \nu_n \tau_n^{-1}= 1-(e+s_1) \tau_n^{-1}$, which yields  \\
 \ind $v_{n+1}+ \dxnp = \gamma_n (v_n+\dxn) + (\ttn-\gamma_n ) \dxn$. \\
 Then, setting $H_n = v_n+ \dxn$, we equivalently have \\
 \ind  $H_{n+1}  = \gamma_n H_{n}  + ( \ttn-\gamma_n) \dxn,$ \\
 along with $\ttn-\gamma_n= - (e-s_0+s_1)\tau_n^{-1}$,   which  amounts to \\
   \ind $H_{n+1}  = \gamma_n H_{n}  + (1-\gamma_n) \frac{\ttn-\gamma_n}{1-\gamma_n} \dxn$.\\
 Hence,  by convexity of the squared norm we infer that  \\
 \ind $\ba \|H_{n+1} \|^2_M  \le  \gamma_n \|H_{n}\|^2_M    + (1-\gamma_n) \lb  \frac{\gamma_n- \ttn}{1-\gamma_n} \rb ^2 \|\dxn\|^2_M \\
  \hskip 1.5cm  =  \gamma_n \|H_{n}\|^2_M    +  \frac{\lb e-s_0+s_1 \rb ^2}{1-\gamma_n}  \tau_n^{-2}  \|\dxn\|^2_M . \ea $ \\
  Hence we obtain \\
 \ind $\|H_{n+1} \|^2_M   = (1-s_0 \tau_n^{-1} ) \|H_{n} \|^2 _M  + s_0^{-1} (e-s_0+s_1)^2 \tau_n^{-1}  \|\dxn\|^2_M $. \\
 Then multiplying this last inequality by $\tau_n^2$ amounts to \\
 \ind $\tau_n^2\|H_{n+1} \|^2_M   \le  (\tau_n^2-s_0 \tau_n ) \|H_{n} \|^2_M   + s_0^{-1} (e-s_0+s_1)^2 \tau_n  \|\dxn\|^2_M $, \\
 while by $\dot{\tau_n}=s_1$ (from the definitions of $\tau_n$ and $\nu_n$) we simply have \\
 \ind  $\tau_n^2 -  \tau_{n-1}^2 \le  s_1 (\tau_n+ \tau_{n-1}) \le  2 s_1 \tau_n$. \\
Combining these last two results amounts to \\
\ind $\tau_n^2\|H_{n+1} \|^2_M   \le \tau_n^2 \|H_{n} \|^2_M  - (s_0 \tau_n- 2 s_1 \tau_n ) \|H_{n} \|^2_M   + s_0^{-1} (e-s_0+s_1)^2 \tau_n  \|\dxn\|^2_M $, \\
 which leads to the desired inequality. \\
 \ind Next  we prove the second part of the lemma. Clearly, the assumptions of Lemma \ref{estimg2} guarantee that   $\sum_n  \tau_n  \|\dxn\|^2_M < \li$ (according to Lemma \ref{estimg0}).
     Consequently, by (\ref{mst}) under the condition $0 \le s_1 < (1/2)s_0$,   we classically deduce that  $\sum_n  \tau_n  \|v_{n}+ \dot{x}_{n} \|^2_M< \li $ (namely  (\ref{dd1})) and that the sequence ($\tau_{n}^2 \|v_{n-1}+ \dot{x}_{n-1} \|^2_M$) is convergent. Thus, there exists $ l_1 \ge 0$ such that $ \lim_{n \to \li} \tau_{n}^2 \|v_{n-1}+ \dot{x}_{n-1} \|^2_M=l_1$, hence, we  also  have   $ \lim_{n \to \li} \tau_{n}^2 \| v_{n}+ \dot{x}_{n} \|^2_M=l_1$ (since $\frac{\tau_{n-1}}{\tau_{n}}  \to 1$ as $n \to \li$). So, noticing that  $\sum_n \tau_n^{-1} =\li$, we are led to  $l_1=0$, which proves (\ref{dd2})\edem \\
 %
  %
 %
 %
\subsection{\label{pthmg} Proof of Theorem \ref{thmg}.}
Let us begin with  observing  that the conditions (\ref{cpal1}) and (\ref{cpal2}) correspond to condition (\ref{cic}) with $i_c=1$ and $i_c=2$, respectively.  The  rest of the proof will be divided into the following two steps {\rm \bf (B1),  (B2)} and {\rm \bf (B3)}: \\
\\
 \ind {\rm \bf (B1)} In order to reach  (\ref{borg1}), (\ref{borg2}) and (\ref{borg1b}), we prove the following estimates:
\begin{eqnarray}
\label{ee}
\label{dd3}
&& \mbox{$\sum_n n   \|v_n \|^2_M < \li$}, \\
\label{dd4}
&& \mbox{$\sum_n   n  \|\la v_n, \dxnp \ra\| _M   < \li$}, \\
\label{dd5}
&& \mbox{$\|\dxn\|_M=o(n^{-1})$}, \\
\label{dd6}
&& \mbox{$\|v_n\|_M=o(n^{-1})$}.
\end{eqnarray}
\ind \ind  Indeed, from a quick computation, we have
 \be  \mbox{$n  \| v_n\|_M^2 \le 2 n  \| v_n+ \dxn\|_M^2 +  2 n  \|\dxn\|_M^2$}. \ee
So, by $\sum_n n \|\dxn\|_M^2< \li$ (from (\ref{pp4})) and  $\sum_n n \| v_n+ \dxn\|_M^2< \li$ (from (\ref{dd1})), we immediately obtain  (\ref{dd3}). The estimate (\ref{dd4}) is an immediate consequence of $\sum_n n \|\dxn\|^2< \li$ (from (\ref{pp4})) and
$\sum_n n \| v_n\|^2_M < \li$ (from (\ref{dd3})).
Next, passing to the limit as $s \to 0^+$
  in (\ref{dju2}) amounts to
 \be \label{dju3}
   \ba \nu_{n+1}^2 \|\dxnp\|^2_M  -  \nu_{n}^2 \|\dxn\|^2_M
    +    s_0 (e+\nu_{n+1})  \la  \vn, \dxnp \ra_M
   +   \frac{1}{2} e   ( e+2  \nu_{n+1} )  \|\dxnp\|^2 _M    \le   0.
    \ea   \ee
   Then, in light of
$\sum_n   n   \|\la v_n, \dxnp \ra\|_M    < \li$ (from (\ref{dd4})) and $\sum_n n   \|\dxn\|^2_M < \li$ (from (\ref{pp4})), we
 derive,   from (\ref{dju3}),  that ($\nu_{n} ^2 \|\dxn\|^2$) is convergent, namely, there exists some $ l_2 \ge 0$ such that  $ \lim_{n \to \li} \nu_{n} ^2 \|\dxn\|^2_M=l_2$. Moreover,  by  $\sum_n  \nu_n \|\dxn\|^2_M < \li$, and recalling that   $\sum_n \nu_n^{-1}=\li$, we get $ \liminf_{n \to \li} \nu_{n} ^2 \|\dxn\|^2_M=0$. It follows that
    $l_2=0$, that is (\ref{dd5}). Next combining this last result with $\|\dxn+ v_n\|_M=o(n^{-1})$ (from (\ref{dd2})) gives us
    $ \lim_{n \to \li} n   \|v_n\|_M=0$, that is (\ref{dd6}). It can be seen that the estimations in   (\ref{borg1}) are  given by  (\ref{dd5}), (\ref{pp4}) and (\ref{pp6}), respectively.
     The estimates in  (\ref{borg2}) follow from (\ref{dd6}) and (\ref{dd3}). Furthermore, by   $\yn=\xn+ \lb 1- \frac{1}{w} \rb \vn$, we obviously have $\dyn=\dxn+ \lb 1- \frac{1}{w} \rb \dot{v}_n$, which implies that
      \be \label{bido1}  \mbox{$\|\dyn\|^2_M \le 2 \|\dxn\|^2_M + 2 \lb 1- \frac{1}{w} \rb^2 \|\dot{v}_n \|^2_M$}. \ee
      Therefore, by (\ref{bido1}) in light of $\|\dxn\|^2=o(n^{-2})$ (from (\ref{dd5})) and   $ \sum_{n }n ^2 \|\dot{v}_n \|^2_{ {M} } < \li$ (from (\ref{pp3})), we obtain  $\|\dyn\|^2=o(n^{-2})$, that is the first result in  (\ref{borg1b}). In addition, by (\ref{bido1}), along with   $ \sum_n n \|\dxn\|^2 < \li$ (from (\ref{dd4})) and $ \sum_{n }n ^2 \|\dot{v}_n \|^2_{ {M} } < \li$ (from (\ref{pp3})), we are led to $ \sum_n n \|\dyn\|^2 < \li$, that is the second  result in  (\ref{borg1b}).\\
     \\
\ind {\rm \bf (B2)} Let us prove  (\ref{borg3a}) and (\ref{borg3c}).      From (\ref{borg2}) and $G_{\ld}^M (z_n) = (\ld w)v_{n+1}$, we readily have
\be \label{cab0} \mbox{$\|G_{\ld}^M(z_n)\|_M^2 ={\cal O}(n^{-2})$ and $\sum_nn  \|G_{\ld}^M(z_n)\|_M^2 <\li$}.\ee
 In addition, given $\{x,z \}  \subset \Hs$ and setting $\Delta G_{\ld}^M(x, z)= G_{\ld}^M(x)- G_{\ld}^M(z)$, by the simple decomposition
       $G_{\ld}^M(x)= \Delta G_{\ld}^M(x, z) + G_{\ld}^M(z)$,
        we readily get
    \be  \label{cab1}  \mbox{$\|G_{\ld}^M(x)\|_M^2 \le  2\|\Delta G_{\ld}^M(x, z) \|^2_M + 2 \|G_{\ld}^M(z) \|_M^2$}. \ee
      Moreover, from condition (\ref{cic}) and (\ref{mohga}) we have
      \be
 \hskip 1.cm     \ld \la  H_{i_c}\Delta G_{\ld}^M(x, z) , \Delta G_{\ld}^M(x, z)  \ra \le \la \Delta G_{\ld}^M(x, z)  , x-z \ra _M  ,
 \ee
 where $  H_{i_c}=M-K_{i_c} I$ (with $K_{i_c}=\delta I $ if $i_c=1$ and $K_{i_c}=\ld L$ otherwise), hence, by applying Peter-Paul's inequality, we obtain \\
  \ind $ \ld \la  H_{i_c}\Delta G_{\ld}^M(x, z) , \Delta G_{\ld}^M(x, z)  \ra \le  (1/2)(\ld w^2)  \| \Delta G_{\ld}^M(x, z)\|^2_M  +  (1/2)  (\ld w^2)^{-1}\|x-z \|_M ^2$, \\
 or equivalently
 \be
 \hskip 1.cm     \ld  \la \lb   H_{i_c} -   (1/2) w^2 M \rb \Delta G_{\ld}^M(x, z) , \Delta G_{\ld}^M(x, z)  \ra \le   (1/2)  (\ld w^2)^{-1} \|x-z \|_M ^2.
 \ee
Furthermore, we simply have \\
\ind $\ba  H_{i_c} -   (1/2) w^2 M= (1-w) M+ (w-w^2) M -K_{i_c} I+  (1/2)w^2 M \\
\hskip 2.cm =(1-w) M+ w(1-w)\bar{M}_{i_c}+  (1/2)w^2M \\
\hskip 2.cm = (1/2)(1+ (1-w)^2) M+ w(1-w)\bar{M}_{i_c} . \ea $  \\
Then, as $\bar{M}_{i_c}$ is positive definite (from condition (\ref{cic})),  combining the previous two results  yields
\be
 \hskip 1.cm   \ld \| \Delta G_{\ld}^M(x, z)\|^2_M  \le  (\ld w)^{-1} \|x-z \|_M ^2,
 \ee
which by (\ref{cab1}) entails that
 \be  \label{cab2}  \mbox{$\|G_{\ld}^M(x)\|_M^2 \le  2  \ld^{-2}   w^{-1} \|x-z\|_M^2 + 2\|G_{\ld}^M(z) \|_M^2 $}. \ee
 In particular,  using (\ref{cab2}), by $\xnp -z_n=-\ld w G_{\ld}^M(z_n)$ we obtain
  \be  \label{cab2a}   \mbox{$\|G_{\ld}^M(\xnp)\|_M^2 \le  2  (w+1)\|G_{\ld}^M(z_n) \|_M^2 $}, \ee
  which  in light of (\ref{cab0}) amounts to (\ref{borg3a}).  Again using  (\ref{cab2}),  by $\ynp=\xnp+ \lb 1- \frac{1}{w} \rb \vnp$, we obtain \\
    \ind $\|G_{\ld}^M(\ynp)\|_M^2 \le  2  \ld^{-2}   w^{-1}\lb 1- \frac{1}{w} \rb^2 \|\vnp\|_M^2 + 2\|G_{\ld}^M(\xnp) \|_M^2 $,\\
 which  in light of  (\ref{borg2}) and (\ref{borg3a}) implies  (\ref{borg3c}).  \\
    \\
 \ind {\rm \bf (B3)} Let us prove (\ref{borg5})  and  (\ref{borg6}).   From Proposition \ref{jeu1}, we know that there exists a sequence  $(y_n^*)\subset \Hs$ (given by (\ref{jeu1})) verifying
 \be \label{bido0} \mbox{$y_{n}^* \in (A+B)(y_{n})$}, \ee
 where  $y_n=x_n+ \lb 1-\frac{1}{w} \rb v_n$. It can also be noticed from  (\ref{piss2}) that $\|y_n^*\|_M = {\cal O}(n^{-1})$, since $\|v_n \|_M ={\cal O}(n^{-1})$ (from (\ref{dd6})). This leads to  (\ref{borg5}). At once, we prove  (\ref{borg6}),  by means of the well-known Opial lemma which guarantees that $(x_n)$ converges to some element of $S$, provided that the following results hold: \\
    \ind \ind (h1) for any $q \in S$, the sequence ($\|x_n-q\|_M$) is convergent, \\
    \ind \ind (h2) any weak-cluster point of ($x_n$), in ($\Hs, \|.\|_M$), belongs to $S$}. \\
 \ind Let us prove (h1).  Take $q \in S$. Clearly, as  a straightforward consequence of the bounded-ness of $(x_n)$ (given by  (\ref{pp0})) along with (\ref{dd6}) we have
\be \label{gris} \la v_n, \xn-q \ra_M =o(n^{-1}). \ee
Moreover, we know that  (${\cal E}_n(s_0,q)$) is convergent (from Lemma \ref{estimg0}) and that it writes
 \be \label{lf3} \es \es  \ba {\cal E}_n(s_0,q)=  \lb \frac{1}{2} \rb   \| s_0  (q-\xn)- \nu_n \dxn  \|^2_M +  \lb \frac{1}{2} \rb \beta_{s_0} \|\xn -q\|^2_M + s_0( e+2  \nu_{n} )  \la v_n, \xn-q \ra_M ,
   \ea \ee
  where $\beta_{s_0}=s_0(e-s_0)$. Then, by  $\nu_n \|\dxn\|_M \to 0$ (from (\ref{dd5})) and $( e+2  \nu_{n} )  \la v_n, \xn-q \ra_M \to 0$ (according to  (\ref{gris})) as $n\to \li$, we deduce that
  \be \label{lf4} \es \es  \ba \lim_{n \to \li}{\cal E}_n(s_0,q)=  \lim_{n \to \li}  \frac{1}{2} s_0 e  \|\xn -q\|_M^2.
   \ea \ee
   This entails (h1). Now, we prove (h2). Let  $u$ be a weak cluster point of $(x_n)$ in ($\Hs, \|.\|_M$), namely  there exists a subsequence ($x_{n_k}$) that converges weakly to $u$ in ($\Hs, \|.\|_M$), as $k \to \li$. Observe that, as $n \to \li$, by $y_n=x_n+ \lb 1-\frac{1}{w} \rb v_n$ and $\|v_n \|_M \to 0 $ (from (\ref{dd6})) we  have  $\|y_n -x_n\|_M \to 0 $ , whence,   ($y_{n_k}$)   converges weakly to $u$ in ($\Hs, \|.\|_M$) (as $k \to \li$).
 Moreover,  by (\ref{borg5}) we know that $\|y_{n}^*\|_M \to 0 $ (as $n \to \li$), while   (\ref{bido0}) gives us
 \be \label{bido2} \mbox{$y_{n_k}^* \in (A+B)(y_{n_k})$}. \ee
  Then passing to the limit as $k \to \li$ in (\ref{bido2})
   and recalling that the graph of the  maximally monotone operator $A+B$ is sequentially closed with respect to the weak-strong topology of the product space $\Hs \mul \Hs$ (see, for instance,  \cite{brez}), we deduce that $0 \in (A+B)(u)$, namely $ u \in S$. This proves (h2) and completes the proof\edem \\
%

\section{\label{pado3} \large From (CRIFBA) to a  generalized variant (G-CRIFBA).}
\setcounter{equation}{0}
Our purpose here is to adapt  CRIFBA to the  general structured monotone inclusion problem
\be \label{gmon} \mbox{find $\bar{x} \in S:= \lb B +  \sum_{k=1}^p A_k \rb^{-1}(0) \neq \emptyset$}, \ee
where $B: \Hs \to \Hs$ is  $\beta$-co-coercive on $\Hs$, while $(A_i)_{i=1}^p: \Hs \to \Hs$ is a family of $p$  maximally  monotone operators  whose resolvent operators  can be easily evaluated. For the sake of clarity we do not include pre-conditioning in the proposed method. To deal with this problem, we follow the methodology of  Raguet-Fadili-Peyre \cite{rafape}, as described  below. \\
\\
\ind  Let   $\{\rho_k \}_{k=1}^p \subset (0,1)$ be  such that $\sum_{k=1}^p \rho_k=1$ and consider the Hilbert space $E=\Hs^p$ endowed with the scalar product $\lb .|. \rb_{\Hs^p}$  defined,  for elements $x=(x_k)_{k=1}^p$ and  $y=(y_k)_{k=1}^p$ belonging to  $E$,  by $\lb x |y \rb_{\Hs^p}= \sum_{k=1}^p \rho_k\la x_k,y_k \ra$.  The induced   norm of $\lb .|. \rb_{\Hs^p}$ will be   denoted by $\| . \|_{\Hs^p}$.\\
\ind Consider also the auxiliary problem
\be \label{gmon2}\mbox{find $\bar{z} \in {S}_p:=\{ z \in \Hs^p \esc ; \esc \sum_{k=1}^p \rho_i z_i \in  S \}$},\ee
  which was shown to have a nonempty solution set $S_p$ (whenever $S \neq \emptyset$). It can also be reformulated as a monotone inclusion   that fits  the structure (\ref{pbg}) and (\ref{cdi}) on $E$.
Introduce indeed  the    mappings  $\bar{A}_G$ and $ \tilde{B}_G$  from $E$ onto  $E$ defined for $(x_i )_{i=1}^p \in E$ by
\be
  \mbox{$  \bar{A}_G \bigg( (x_i )_{i=1}^p \bigg) =  \bigg( \frac{\ld }{\rho_i} A_i(x_i) \bigg)_{i=1}^p, $} \es \es
 \mbox{ $ \tilde{B}_G \bigg( (x_i ) _{i=1}^p \bigg)=  \bigg( B(x_i) \bigg)_{i=1}^p,$}
\ee
 and let  $N_{\Gamma}: \Hs^p \to 2 ^{\Hs^p}$ be  the normal cone to the (nonempty) closed convex set \\
  \ind $\Gamma=\{ \big(x_i\big)_{i=1}^p \in \Hs^p \esc ;  \esc x_1=x_2=...=x_p \}$. \\
  It can be checked that  $\bar{A}_G$ and $\tilde{B}_G$ are maximally monotone operator on $E$. So the
  reflection operators $R_{\bar{A}_G}= 2 J_{\bar{A}_G }-I_{\Hs^p}$ and  $R_{ N_{\Gamma} }= 2 J_{N_{\Gamma} } -I_{\Hs^p}$ are well-defined, which   allows us to consider the mappings   $T_{1}$,  $T_{2}$ and $T$, from $E$ onto $E$ and   such that
   \be  \label{gmono} \mbox{ $ T_{1} = \frac{1}{2} (R_{ \bar{A}_G} \circ R_{N_{\Gamma}} +I_{\Hs^p})$, 
\es $ T_{2} = I_{\Hs^p} - \ld  \tilde{B}_G\circ  J_{N_{\Gamma}}$ \es and \ind  $T=T_{1} \circ T_{2}$}. \ee
It is established in  \cite{rafape} the results given in the next proposition. \\
\begin{proposition} The following statements are obtained: 
\begin{itemize}
\item   (See \cite[Propositions 4.1, 4.2 and 4.6]{rafape}) There exists some  maximally monotone operator ${\cal A}: E \to 2 ^{E}$ such that
{\small \be  \label{gmon4} \ba \mbox{$S_p= \lb {\cal A}+ \tilde{B}_G \circ J_{N_{\Gamma}} \rb ^{-1}(0)$}. \ea \ee}
\item (See  the proof of \cite[Proposition 4.1]{rafape}) The operator  $\tilde{B}_G \circ J_{N_{\Gamma} }$ is   $\beta$ co-coercive. 
\item   $S_p$ is nothing but the fixed point set of the operator $T$ which can be rewritten as 
{\small\be  \label{gmon3} \ba T= J_{{\cal A}}\circ ( I_{\Hs^p} - \ld  \tilde{B}_G \circ  J_{N_{\Gamma}}). \ea \ee}
\end{itemize}
\end{proposition}
 Consequently, a strategy to solve (\ref{gmon}) consists first of approaching  an element of $S_p$  (namely a fixed point of $T$)  by means of  sequences $(z_n)=\bigg( (z_{n,k})_{k=1}^p \bigg)  \subset \Hs^p$ and $(\zeta_n)=\bigg( (\zeta_{n,k})_{k=1}^p \bigg)  \subset \Hs^p$ generated by CRIFBA, in the context of    (\ref{gmon4})-(\ref{gmon3}), as follows:
 \begin{subeqnarray}
&  & \hskip -1.cm \mbox{$z_{n}= \zeta_{n} + \ttn   (\zeta_{n} - \zeta_{n-1}) +  \gamma_n  (z_{n-1}- \zeta_{n}) $},\\
& &  \hskip -1.cm \mbox{$\zeta_{n+1}=(1-w) z_{n} + w T(z_n)$}.
  \end{subeqnarray}
 Next, we  derive  an element of $S$ as  the limit of  $(x_n) \subset \Hs$ given by
 $x_n=\sum_{k=1}^p \rho_k \zeta_{n,k}$. This leads us (see the proof of Theorem \ref{saba1})  to  the  algorithm (G-CRIFBA) given below  : \\
\\
{\bf (G-CRIFBA)}: \\
 $\rhd$ {\bf Step 1} (\un{initialization}): \\
 \ind \ind Let $\{ z_{-1}, \zeta_{-1}, \zeta_{0} \} \subset \Hs^p$,   \es   $\{e, s_0, s_1,   \nu_0, \ld, w \} \subset [0,\li)$,
 \es $\{\rho_k\}_{k=1}^p \subset (0,1)$,\\
\ind \ind and  set $\nu_n=s_1n+\nu_0$, \es $\ttn= 1- \frac{e+s_1 }{e+ \nu_{n+1}}$ \es and  \es $\gamma _n=  1-\frac{s_0}{e+ \nu_{n+1}} $. \\
 $\rhd$ {\bf Step 2} (\un{main step}): \\
\ind \ind Given  $\{ z_{n-1}, \zeta_{n-1}, \zeta_{n} \} \subset \Hs^p$  (with $n \ge 0$), we   compute
 \begin{subeqnarray}
 \label{filo}
  \label{filo1}
&  & \hskip -0.cm \mbox{$z_{n}= \zeta_{n} +  \ttn   (\zeta_{n} - \zeta_{n-1}) +  \gamma_n  (z_{n-1}- \zeta_{n})  $},\\
 \label{filo2}
&  & \hskip -0.cm \mbox{$u_{n}=\sum_{k=1}^p \rho_k z_{n,k}  $},\\
 \label{filo3}
 & &  \hskip -0.cm \zeta_{n+1}= \bigg(  z_{n,k} + w \lb  J _{\frac{\ld }{\rho_k}A_k} \lb 2  u_{n} -\ld  B(u_{n})- z_{n,k} \rb -   u_{n} \rb  \bigg) _{k=1}^p, \\
 \label{filo4}
& &  \hskip -0.cm \mbox{$x_{n+1}=\sum_{k=1}^p \rho_k \zeta_{n,k}$}.
  \end{subeqnarray}
 \ind The next  theorem sets the convergence rates of  the iterates $(\zeta_n)$ generated by G-CRIFBA in terms of  discrete velocity and fixed point residual to $\| \dot{\zeta}_n \|_{\Hs^p} =o( n^{-1})$ and $ \|T({\zeta}_n ) -{\zeta}_n \|_{\Hs^p} =o( n^{-1})$, respectively, instead of the  rates $\| \dot{\zeta}_n \|_{\Hs^p} ={\cal O}( n^{-1/2})$ and $ \|T({\zeta}_n ) -{\zeta}_n \|_{\Hs^p} = {\cal O}( n^{-1/2})$ obtained for classical fixed point iterations of $T$ as in \cite{rafape} (that is   (\ref{filo}) with $w=1$ and $\theta_n=\gamma_n=0$). \\
 \begin{theorem} \label{saba1}
  Let
  $\{x_n\} \subset \Hs $ be   generated by G-CRIFBA  with      $\{\rho_k \}_{k=1}^p \subset (0,1)$ verifying $\sum_{k=1}^p \rho_k=1$, together with the other parameters such that 
   \begin{eqnarray}
   && \mbox{$0 < \ld  <  4w(1-w)\beta$, \es $0 < w < 1$, \es $2s_1 <   s_0 <e$}.
   \end{eqnarray}
Then the  following  results  are reached:
   \begin{subeqnarray}
   \label{gorg}
    \label{gorg1}
   &  &     \mbox{$\| \dot{\zeta}_n \|_{\Hs^p} =o( n^{-2})$,
   \es     $\sum_n n  \| \dot{\zeta}_n  \|^2_{\Hs^p} <  \li$, \es $\sum_n n^2  \| \dot{\zeta}_{n+1} -\dot{\zeta}_n  \|^2_{\Hs^p} <  \li$}, \\
   \label{gorg2}
   &  &    \mbox{$ \|\zeta_{n+1}-z_n \|^2 _{\Hs^p}=o( n^{-2})$,
   \es $\sum_n n   \|\zeta_{n+1}-z_n\|_{\Hs^p} ^2 <  \li$}, \\
    \label{gorg3}
   &  &    \mbox{$ \|T(\zeta_n) -\zeta_n \|_{\Hs^p} =o( n^{-1})$,
   \es $\sum_n n   \|T(\zeta_n) -\zeta_n \|_{\Hs^p} ^2 <  \li$}, \\
  \label{gorg4}
 &  &  \mbox{$\exists \bar{\zeta}\in S_p$, s.t. (for $k=1,..,p$) $\bar{\zeta}_{n,k}  \rightharpoondown \bar{\zeta}_k$  weakly in $\Hs$, as $n \to \li$,}\\
   \label{gorg5}
   &  &  \mbox{$x_n  \rightharpoondown \bar{x}=\sum_{k=1}^p \rho_k  \bar{\zeta}_{k} \in S$  weakly in $\Hs$}.
   \end{subeqnarray}
\end{theorem}
\ddem
Let us evaluate the operator  $T$ on $E$  from its formulation given by (\ref{gmono}), namely $T=T_{1} \circ T_{2}$. It  can be  checked (see \cite[Lemma 4.1]{rafape}) that,  for    $( x_k) _{k=1}^p \in E$, we have
\begin{eqnarray}
&&   \mbox{$ J_{N_{\Gamma}}\big(\{x_i \}_{i=1}^p \big)=   \big( \sum_{j=1}^p \rho_j x_j \big)_{i=1}^p, $} \\
&& \mbox{$R_{ N_{\Gamma} }\big(\{ x_k\}_{k=1}^p \big)=  \big( 2 \sum_{i=1}^p \rho_i x_i - x_k \big) _{k=1}^p$}, \\
&& \mbox{$R_{\bar{A}_G  } \big(\{x_k\}_{k=1}^p \big)=   \big( 2 J_{\frac{\ld }{\rho_k} A_k}(x_k) - x_k \big) _{k=1}^p $}.
\end{eqnarray}
Then,  for    $\big( z_k\big)_{k=1}^p \in E$, by $ T_{2 } = I_{\Hs^p} - \ld  \tilde{B}_G\circ  J_{N_{\Gamma}}$ we obtain  \\
\ind  $\ba T_{2 }\big( ( z_k) _{k=1}^p \big)=  ( z_k)_{k=1}^p - \ld  ( \tilde{B}_G \circ J_{N_{\Gamma}})\big( ( z_k )  _{k=1}^p \big) \\
\hskip 1.5cm = ( z_k) _{k=1}^p -  \ld  \tilde{B}_G \bigg( \big( \sum_{i=1}^p \rho_iz_i  \big) _{k=1}^p \esc \bigg) \\
\hskip 1.5cm = \bigg( z_k-  \ld  \esc {B}\big(  \sum_{i=1}^p \rho_iz_i  \big) \esc \bigg)  _{k=1}^p.\\
\ea $\\
In addition, for  $(y_k) _{k=1}^p \in E$ and setting $\bar{y}=\sum_{i=1}^p \rho_iy_i$, by \\
 \ind $ T_{1} = (1/2) (R_{ \bar{A}_G} \circ R_{N_{\Gamma}} +I_{\Hs^p})$ \\
we get  \\
\ind $\ba T_{1}\big( ( y_k)_{k=1}^p \big)=\frac{1}{2} \bigg( R_{ \bar{A}_G}\lb   \lb  2 \bar{y} -y_k\rb _{k=1}^p \rb  + ( y_k)_{k=1}^p \bigg)  \\
\hskip 1.5cm =\frac{1}{2} \bigg(  \lb  2 J _{\frac{\ld }{\rho_k}A_k}\lb 2 \bar{y}  -y_k \rb - 2 \bar{y}  +y_k  \rb  _{k=1}^p  + ( y_k)_{k=1}^p \bigg)  \\
\hskip 1.5cm =  \bigg(  J _{\frac{\ld }{\rho_k}A_k}\lb 2  \bar{y} -y_k \rb -  \bar{y} + y_k  \bigg) _{k=1}^p .  \\
\ea $\\
Hence, taking  $ ( y_k)_{k=1}^p := T_{2 }\big( ( z_k) _{k=1}^p \big)$ and setting $\bar{z}=\sum_{i=1}^p \rho_iz_i$, we have\\
 \ind $ T( ( z_k)_{k=1}^p )   =   T_{1 } \lb ( y_k)_{k=1}^p \rb$, \\
 or equivalently    \\
 \ind $ \ba
 T( ( z_k)_{k=1}^p )  =  T_{1 } \bigg( \big(   z_{k} -  \ld  {B}\lb \bar{z}\rb \big) _{k=1}^p \bigg)  = \bigg(   J _{\frac{\ld }{\rho_k}A_k} \lb 2 \bar{z}- \ld  {B} \lb \bar{z}
   \rb -z_{k} \rb -  \bar{z} + z_{k} \bigg) _{k=1}^p. \ea $ \\
   It follows that
 \be
 (1-w) ( z_k)_{k=1}^p+ w T\big( ( z_k)_{k=1}^p \big)
  = \bigg(  z_k + w  \lb   J _{\frac{\ld }{\rho_k}A_k} \lb 2 \bar{z}- \ld  {B} \lb \bar{z}
   \rb -z_{k} \rb -   \bar{z} \rb  \bigg) _{k=1}^p. \ee
This leads us to the formulation of  G-CRIFBA and the  results (\ref{gorg1}) to (\ref{gorg4}) follow straightforwardly  from  Theorem \ref{thmg} (also see Remark \ref{remg}), while  (\ref{gorg5}) is immediately deduced  from (\ref{gorg4})\edem \\
 %
%
\section{\label{pado4} \large An application of (CRIFBA) to some convex-concave saddle-point problem.}
\setcounter{equation}{0}
 In this section we apply CRIFBA to the   problem below discussed  by Lorenz-Pock \cite{lp}.\\
   \ind Let $X$  and $Y$  be two Hilbert spaces endowed with scalar products $\la  .,.\ra _X$ and $\la .|.\ra _Y$, respectively,  and induced norms denoted by  $\|.\|_X$ and $\|.\|_Y$, and consider  the following
 saddle-point problem
\be \label{gko} \min_{x \in X} \max_{y \in Y} \es G(x)+ Q(x)+ \la  Kx , y \ra _Y-F^*(y)-P^*(y), \ee
where   $K:X \to Y$ is linear and bounded,   $G: X \to (-\li,\li]$ and  $F^*:Y \to   (-\li,\li]$ are convex functions, while 
 $Q: X \to (-\li,\li]$ and  $P^*:Y \to   (-\li,\li]$ are convex differentiable functions with Lipschitz continuous gradient (whose respective Lipschitz constants are $l_Q$ and $l_{P^*}$). \\
 \\
 \ind The above problem covers several primal-dual formulation of  nonlinear  problems encountered  for instance in  image processing.\\
 \\
 \ind We denote by $S$ the solution set of (\ref{gko}) and we assume that
 $S\neq \emptyset$. \\
 \\
 \ind Introduce  the Hilbert space $E=X \mul Y$ endowed with the scalar product $( .|. )$ defined for  $ \zeta_1=(x_1,y_1) \in E$ and $\zeta_2=(x_2,y_2)  \in  E$ by  $( \zeta_1 | \zeta_2 )= \la x_1, x_2 \ra _X+\la y_1, y_2 \ra _Y $,
and let us denote its induced norm by $\|.\|$.
%
%
So,   (\ref{gko}) through its optimality condition   can be re-formulated  as
\be \label{pvsp} \mbox{find $(x,y) \in E$ such that $0 \in (A+B) \left ( \begin{array}{c}
x   \\
y      \\
\end{array} \right )$}, \ee
 where  $A$ and  $ B$ are the monotone operators on $E$ defined by
{\small \be \label{vilo1} A= \left ( \begin{array}{cc}
\der G   &  K^* \\
-K     & \der F^* \\
\end{array} \right ), \es B= \left ( \begin{array}{cc}
\nabla Q   &  0 \\
0     & \nabla P^* \\
\end{array} \right ). \ee}
It is    established in \cite{lp} the following result. 
\begin{proposition}
(See \cite[Proof of Theorem 5]{lp}) The operator $B$ is co-coercive w.r.t. to the mapping
{\small\be L= \left ( \begin{array}{cc}
l_Q I_X   &  0 \\
0     & l_{P^*} I_Y \\
\end{array} \right ), \ee}
where $I_X$ and $I_Y$ denote the identity mappings on $X$ and $Y$, respectively. \\ \\
\end{proposition}
\ddem Given $\{ \zeta_1=(x_1,y_1), \zeta_2=(x_2,y_2) \} \subset E$, we have \\
 \ind $ \ba \lb  B(\zeta_1)- B(\zeta_2) | \zeta_1-\zeta_2 \rb  \\
\hskip 1.cm =  \la  \nabla Q (x_1)-\nabla Q (x_2) , x_1-x_2 \ra_X +  \la  \nabla P^* (y_1)-\nabla P^*(y_2) | y_1-y_2 \ra_Y \\
\hskip 1.cm\ge l_{Q}^{-1} \|Q (x_1)-\nabla Q (x_2)\|_X + l_{P^*}^{-1} \|\nabla P^* (y_1)-\nabla P^*(y_2)\|_Y \\
\hskip 1.cm = \lb  L^{-1} (B(\zeta_1)- B(\zeta_2)) | B(\zeta_1)- B(\zeta_2)\rb = \| B(\zeta_1)- B(\zeta_2)\|^2_{L^{-1}}\mbox{\edem} \ea$ \\
\\
\ind As a consequence, the  above  monotone inclusion (\ref{pvsp}) enters  the setting of (\ref{pbg}) and (\ref{cdi}) and so it can be solved by means of  the proposed method CRIFBA.
In general, as explained in  \cite{lp}, evaluating the proximal mapping $(I + \ld A)^{-1}$ may be prohibitively expensive, which would make  our algorithm impracticable in the standard case when  $M=I$.  Fortunately, this drawback  can be overcame when choosing the pre-conditioner  mapping $M$ in order   to cancel out the upper off-diagonal block in the  sum $M+A$, as follows
\be \label{vilo2} M= \left ( \begin{array}{cc}
 \tau^{-1} I_X   &  -K^* \\
-K    & \sigma^{-1} I_{Y} \\
\end{array} \right ), \es   M+A= \left ( \begin{array}{cc}
 \tau^{-1} I_X+ \der G   &  0 \\
-2 K    & \sigma^{-1} I_Y+ \der F^* \\
\end{array} \right ),    \ee
where $\tau$ and $\sigma$ are positive real numbers. Furthermore, an easy computation gives us the following result. \\
\begin{proposition}
For any $(\xi',\chi') \in X \mul Y$, we obtain
\be \label{vilo3} (M+A)^{-1}  \left ( \begin{array}{c}
\xi'   \\
\chi '       \\
\end{array} \right )=  \left ( \begin{array}{l}
{\rm prox}_{\tau G}( \tau \xi' )   \\
{\rm prox}_{\sigma F^*} \lb \sigma \chi ' + 2 \sigma K {\rm prox}_{\tau G}( \tau \xi' ) \rb      \\
\end{array} \right ).
 \ee
 \end{proposition}
\ind Clearly,  a  solution to  (\ref{pvsp}) can be approximated    by means of a sequence $\{(x_n,y_n)\}  \subset X \mul Y$ generated by CRIFBA  (with $\ld=1$), in the context of  (\ref{vilo1}) and  (\ref{vilo2}), as follows
\begin{subeqnarray}
&  & \hskip -0.cm \mbox{$(\xi_n, \chi_n)= (x_{n}, y_{n}) +  \ttn   ((x_{n}, y_{n})- (x_{n-1}, y_{n-1}))
+  \gamma_n ((\xi_{n-1}, \chi_{n-1})-(x_{n}, y_{n})  )   $},\\
& &  \hskip -0.cm \mbox{$(x_{n+1},y_{n+1}) =(1-w)(\xi_n, \chi_n)+ w T (\xi_n, \chi_n)   $},
  \end{subeqnarray}
where $T = (M+A)^{-1}\lb M-B   \rb$.  This in light of (\ref{vilo3})  leads us to the following corrected relaxed inertial primal dual  algorithm: \\
\\
{\bf (CRIPDA)}: \\
 $\rhd$ {\bf Step 1} (\un{initialization}): \\
 \ind \ind Let $\{ (\xi_{-1},\chi_{-1}), (x_{-1},y_{-1}), (x_{0},y_{0}) \} \subset X\mul Y$, \ind     $\{e, s_0, s_1,   \nu_0,  w, \sigma, \tau \} \subset [0,\li)$, \\
 \ind \ind and set $\nu_n=s_1n+\nu_0$, \es $\ttn= 1- \frac{e+s_1 }{e+ \nu_{n+1}}$ \es and
 \es  $\gamma _n=  1-\frac{s_0}{e+ \nu_{n+1}}$. \\
 $\rhd$ {\bf Step 2} (\un{main step}): \\
\ind \ind Given  $\{ (\xi_{n-1}, \chi_{n-1}), (x_{n-1},y_{n-1}),  (x_{n}, y_{n}) \} \subset X \mul Y$ (with $n \ge 0$), we   compute
 \begin{subeqnarray}
&  & \hskip -1.cm \mbox{$\xi_{n}= x_{n} + \ttn   (x_{n} - x_{n-1}) +  \gamma _n  (\xi_{n-1}-x_n) $},\\
&  & \hskip -1.cm \mbox{$\chi_{n}= y_{n} +  \ttn   (y_{n} - y_{n-1}) +  \gamma _n  (\chi_{n-1}-y_n) $},\\
& &  \hskip -1.cm \mbox{$x_{n+1}=(1-w) \xi_{n} + (w) {\rm prox}_{\tau G} \bigg( \xi_n -\tau (\nabla Q( \xi_n)+  K^* \chi_n) \bigg)$}, \\
& &  \hskip -1.cm \mbox{$\bar{\xi }_{n}=2 w^{-1} \lb x_{n+1}- (1-w) \xi_{n}\rb $}, \\
& &  \hskip -1.cm \mbox{$y_{n+1}=(1-w) \chi_n   + (w) {\rm prox}_{\sigma F^*} \bigg(   \chi_n  -\sigma(  \nabla P^*(\chi_n) -  K  \bar{\xi }_{n}) \bigg)$}.
  \end{subeqnarray}
  \ind  In the absence of any   correction term and  relaxation factor (that is $\gamma_n=0$ and $w=1$) we retrieve the primal-dual algorithm  in \cite{lp}, which was discussed with step-size rules regarding  the momentum term $\ttn$.
   In the absence of inertial and correction terms (that is $\theta_n=\gamma_n=0$) we retrieve the primal-dual algorithms,   proposed by  Condat \cite{condat} (for  $P^*=0$) and Vu \cite{vu}, which were  also investigated with varying relaxation factors. Compared with  these methods, a fast convergence rate is proved for   CRIPDA.  \\
  \\
\ind The next result establishes the convergence of the above algorithm. \\
 \begin{theorem}
  Let
  $\{(x_n,y_n)\} \subset X \mul Y $ be generated by CRIPDA  with  parameters verifying
   \begin{eqnarray}
   && \mbox{$0 < w < 1$, \es $2s_1 <   s_0 <e$}.
   \end{eqnarray}
   Suppose in addition that one of the two  conditions  (\ref{depo1}) and   (\ref{depo2}) holds:
   \begin{eqnarray}
   \label{depo1}
    && \es \es \mbox{\scriptsize$ \delta > \frac{1}{4} \max \{ l_Q,l_{P^*} \} $},\es   \mbox{\scriptsize$\{ \tau, \sigma \} \subset  \lb 0,  \frac{ w(1-w)}{\delta} \rb $},\es   \mbox{\scriptsize $\|K\|^2 < \lb \tau^{-1} -\frac{\delta}{w(1-w)} \rb \lb \sigma^{-1} -\frac{\delta}{w(1-w)}\rb$}, \\
   \label{depo2}
    && \es \es \mbox{\scriptsize$0 <  \tau <   \frac{w(1-w)}{l_Q}$}, \es \mbox{\scriptsize $0 <  \sigma  < \frac{w(1-w)}{l_ {P^*}}$}, \es   \mbox{\scriptsize $\|K\|^2 < \lb \tau^{-1} -\frac{l_{Q}}{w(1-w)} \rb \lb \sigma^{-1} -\frac{l_{P^*}}{w(1-w)}\rb$}.
   \end{eqnarray}
Then the  following  results  are reached:
   \begin{subeqnarray}
   &  &   \es \es \es   \mbox{ $\| (\dxn,\dyn)\|_M^2 =o( n^{-2})$,
       $\sum_n n \| (\dxn,\dyn) \|_M^2 <  \li$, $\sum_n n^2 \| (\dxnp-\dxn,\dynp-\dyn) \|_M^2 <  \li$}, \\
    &  &   \es \es \es   \mbox{ $ \|T(x_n,y_n) -(x_n,y_n) \|_M^2 =o( n^{-2})$,
   \es $\sum_n n  \|T(x_n,y_n) -(x_n,y_n) \|_M^2  <  \li$}, \\
   &  &   \es \es \es \mbox{\es $\exists (\bar{x},\bar{y})  \in S$, s.t. $(x_n,y_n)  \rightharpoondown (\bar{x},\bar{y})  $ weakly in $X \mul Y$}.
   \end{subeqnarray}
\end{theorem}
\ddem
 Clearly, we have  $\|L\|=\sup\{l_{P^*}, l_Q\}$, while it can be checked that the operator $M-\frac{\delta}{w(1-w)}I$ is positive definite if the last two conditions in (\ref{depo1}) are  fulfilled. Moreover, it can be verified that the operator $M-\frac{1}{w(1-w)}L$ is positive definite if (\ref{depo2}) is fulfilled. Therefore, conditions (\ref{cpal1}) and   (\ref{cpal2}) of Theorem \ref{thmg} are satisfied under   conditions (\ref{depo1}) and   (\ref{depo2}), respectively. The rest of the proof is a direct consequence of  Theorem \ref{thmg}\edem \\
%
\setcounter{section}{0}
\renewcommand{\thesection}{\Alph{section}}
\renewcommand{\thesubsection}{\Alph{section} - \arabic{subsection}}
\section{APPENDIX.}
\subsection{\label{pgcd} Proof of Proposition \ref{gcd}.} For simplification reasons, we write $G$ instead of $G^M_{\ld}$ and,   given any mapping $\Gamma:\Hs \to \Hs$ and any elements $\{x_1,x_2\} \subset \Hs$, we denote \\
\ind $\Delta \Gamma(x_1,x_2)=\Gamma(x_1)-\Gamma(x_2)$.\\
Let $\bar{A}=M^{-1}A$,    $\bar{B}=M^{-1}B$ and $C=I-\ld  \bar{B}$. Clearly, $\bar{A}$ is monotone in $(\Hs,|.|_M)$. It is also  obviously seen  for $ x \in \Hs$  that
    $G(x) =\ld^{-1} \bigg( x -J_{\ld \bar{A}}\big(C(x) \big)\bigg) $, hence  \\
    \ind $G(x) =\ld^{-1} \big( x -C(x) \big) + \ld^{-1} \bigg( C(x) -J_{\ld \bar{A}}\big(C(x) \big) \bigg) $, \\
    or equivalently
    \be \label{fiko} \mbox{$G(x) = \bar{B}(x)+ \bar{A}_{\ld}(C(x))$}, \ee
  where $\bar{A}_{\ld}:=\ld^{-1} \lb I -J_{\ld \bar{A}}\rb $ is the Yosida regularization of $\bar{A}$.  Now, given any $(x_1, x_2) \in \Hs^2$, we get \\
   \ind $\ba  \la  \Delta G(x_1,x_2) , x_1-x_2 \ra _M = \la \bar{B}(x_1)+ \bar{A}_{\ld}(C(x_1))-\bar{B}(x_2)- \bar{A}_{\ld}(C(x_2)) , x_1-x_2 \ra _M  \\
   \hskip 3.5cm=   \la \Delta \bar{B}(x_1, x_2) , x_1-x_2 \ra _M  +   \la \Delta (\bar{A}_{\ld} \circ C)(x_1, x_2) , x_1-x_2 \ra _M \\
 \hskip 3.5cm=   \la \Delta \bar{B}(x_1, x_2) , x_1-x_2 \ra _M +   \la \Delta (\bar{A}_{\ld} \circ C)(x_1, x_2) , \Delta C(x_1, x_2) \ra _M \\
 \hskip 4.5cm+  \la \Delta (\bar{A}_{\ld} \circ C)(x_1, x_2)  , \Delta (I-C) (x_1, x_2)  \ra _M, \ea  $\\
 hence, by $I-C=\ld  \bar{B}$, we equivalently obtain
  \be     \label{fiko1} \ba \la \Delta G(x_1,x_2), x_1-x_2 \ra _M  \ge    \la  \Delta \bar{B}(x_1,x_2) , x_1-x_2 \ra _M   +   \la  \Delta (\bar{A}_{\ld} \circ C)(x_1, x_2) , \Delta C(x_1, x_2)  \ra _M \\
 \hskip 4.cm + \ld  \la  \Delta (\bar{A}_{\ld} \circ C)(x_1, x_2) , \Delta \bar{B}(x_1,x_2)  \ra _M.\ea \ee
   Let us estimate separately the last two terms in the right side of the previous inequality.  As a classical result, by the $\ld$-co-coerciveness of $\bar{A}_{\ld}$   in $(\Hs,\|.\|_M)$, we have \\
    \ind   $ \la  \Delta (\bar{A}_{\ld} \circ C)(x_1, x_2) , \Delta {C}(x_1,x_2) \ra _M \ge \ld  \|\Delta (\bar{A}_{\ld} \circ C)(x_1, x_2) \|_M^2$,  \\
    which  by $ \bar{A}_{\ld}\circ C =G- \bar{B}$ (from (\ref{fiko})) can be rewritten as \\
   \ind   $ \la  \Delta (\bar{A}_{\ld} \circ C)(x_1, x_2) , \Delta {C}(x_1,x_2)\ra _M \\
   \ind \ind \ind \ind \ge \ld  \|\Delta G(x_1,x_2)-  \Delta \bar{B}(x_1,x_2) \|_M^2 \\
    \ind \ind \ind \ind =\ld  \|\Delta G(x_1,x_2) |_M^2 + \ld  \|\Delta \bar{B}(x_1,x_2)   \|_M^2
    - 2 \ld \la  \Delta G(x_1,x_2)  ,\Delta \bar{B}(x_1,x_2)  \ra_M  $.  \\
    Moreover,  by $ \bar{A}_{\ld}\circ C =G- \bar{B}$ (from (\ref{fiko})), we simply get \\
  \ind   $ \ba \ld \la  \Delta (\bar{A}_{\ld} \circ C)(x_1, x_2)  , \Delta \bar{B}(x_1,x_2)  \ra _M \\
  \hskip 1.cm = \ld  \la  \Delta G(x_1,x_2) - ( \Delta \bar{B}(x_1,x_2) ) , \Delta \bar{B}(x_1,x_2) \ra _M   \\
 \hskip 1.cm = \ld  \la \Delta {G}(x_1,x_2) ) ,  \Delta \bar{B}(x_1,x_2) \ra _M - \ld \|\Delta \bar{B}(x_1,x_2)  \|^2 _M. \ea  $ \\
    Thus, by (\ref{fiko1}) and the previous arguments, we obtain \\
 \ind $\ba \la \Delta {G}(x_1,x_2), x_1-x_2 \ra _M \\
  \hskip 0.5cm \ge    \la  \Delta \bar{B}(x_1,x_2) , x_1-x_2 \ra _M  +  \ld  \|\Delta {G}(x_1,x_2) \|_M^2 + \ld  \|\Delta \bar{B}(x_1,x_2) \|_M^2 \\
 \hskip 1.cm  - 2 \ld \la \Delta {G}(x_1,x_2),  \Delta \bar{B}(x_1,x_2) \ra_M  + \ld  \la  \Delta {G}(x_1,x_2),   \Delta \bar{B}(x_1,x_2) \ra _M - \ld \| \Delta \bar{B}(x_1,x_2)\|^2 _M \\
\hskip 0.5cm =  \la  \Delta \bar{B}(x_1,x_2) , x_1-x_2 \ra _M  +  \ld  \|\Delta {G}(x_1,x_2) \|_M^2
    -  \ld \la \Delta {G}(x_1,x_2),  \Delta \bar{B}(x_1,x_2) \ra_M.  \\
   \ea$ \\
Hence, reminding that  $\bar{B}=M^{-1}B$, we equivalently obtain  \\
\ind $\ba \la \Delta {G}(x_1,x_2), x_1-x_2 \ra _M   \\
  \hskip 0.5cm \ge     \la \Delta \bar{B}(x_1,x_2) , x_1-x_2 \ra   +  \ld  \|\Delta {G}(x_1,x_2) \|_M^2
    -  \ld \la \Delta {G}(x_1,x_2), \Delta {B}(x_1,x_2) \ra.  \\
  \ea$ \\
  Then by the co-coercivity assumption on $B$ we infer that
  \be \label{moh} \hskip 0.8cm  \la \Delta {G}(x_1,x_2), x_1-x_2 \ra _M \ge     \|\Delta {B}(x_1,x_2) \|^2_{L^{-1}}   +  \ld  \|\Delta {G}(x_1,x_2)  \|_M^2
    -  \ld \la \Delta {G}(x_1,x_2) , \Delta {B}(x_1,x_2) \ra,
 \ee
 that is (\ref{mohg}). \\
 \ind Now, let us prove (\ref{mohg}) for $i=1$. From an easy computation, we  obtain \\
   \ind $ \|\Delta {B}(x_1,x_2)  \|^2 \le    \|L ^{\frac{1}{2}} \|^{2} \|\Delta {B}(x_1,x_2)  \|^2_{L^{-1}}$. \\
  hence,   for any $\delta >0$,  using successively   Peter-Paul's inequality and the previous inequality gives us
  \be  \ba  \la \Delta {G}(x_1,x_2) , \Delta {B}(x_1,x_2) \ra \le  \mbox{$ \delta \|\Delta {G}(x_1,x_2) \|^2$}+  \mbox{$\frac{1}{4 \delta} \|\Delta {B}(x_1,x_2) \|^2$} , \\
   \hskip 4.cm \le
         \mbox{$   \delta \|\Delta {G}(x_1,x_2)\|^2$} +  \mbox{$\frac{1}{4 \delta}  \|L ^{\frac{1}{2}}\|^{2} \|\Delta {B}(x_1,x_2)\|^2_{L^{-1}} $}.
 \ea \ee
  Therefore, combining this last inequality with (\ref{moh})  entails \\
   \ind $\ba \la \Delta {G}(x_1,x_2) , x_1-x_2 \ra _M  \\
  \hskip 1.cm \ge   \lb 1 -  \frac{\ld}{4 \delta }  \|L ^{\frac{1}{2}}\|^{2}  \rb   \|\Delta {B}(x_1,x_2)|^2_{L^{-1}}
     +   \ld  \lb \|\Delta {G}(x_1,x_2) \|_M^2 -   \delta \|\Delta {G}(x_1,x_2)  \|^2\rb \\
    \hskip 1.cm =   \lb 1 -  \frac{\ld}{4 \delta }  \|L ^{\frac{1}{2}}\|^{2}  \rb   \|\Delta {B}(x_1,x_2)\|^2_{L^{-1}}
     +   \ld \la \lb M -\delta I \rb  \Delta {G}(x_1,x_2) , \Delta {G}(x_1,x_2) \ra, \\
  \ea$ \\
  that is  (\ref{mohg}) with $i=1$.  \\
  \ind Let us prove (\ref{mohg}) for $i=2$. Using again  Peter-Paul's inequality, we readily have  \\
  \ind $\ba  \ld \la \Delta {G}(x_1,x_2) , \Delta {B}(x_1,x_2) \ra=\la  \ld \Delta {G}(x_1,x_2) , L ^{\frac{1}{2}} L ^{-\frac{1}{2}}\Delta {B}(x_1,x_2) \ra \\
   \hskip 3.cm=\la  \ld  L ^{\frac{1}{2}}\Delta {G}(x_1,x_2) , L ^{-\frac{1}{2}}\Delta {B}(x_1,x_2) \ra \\
   \hskip 3.cm  \le    \ld^2 \|L ^{\frac{1}{2}}\Delta {G}(x_1,x_2) \|^2+ \frac{1}{4 } \|L ^{-\frac{1}{2}} \Delta {B}(x_1,x_2) \|^2 \\
  \hskip 3.cm  =   \ld^2\la  L \Delta {G}(x_1,x_2), \Delta {G}(x_1,x_2) \ra + \frac{1}{4} \|\Delta {B}(x_1,x_2) \|^2_{L ^{-1} }, \ea $ \\
which, in  light (\ref{moh}),   gives us  \\
 \ind $\ba \la \Delta {G}(x_1,x_2) , x_1-x_2 \ra _M  \ge   \frac{3}{4 }    \|\Delta {B}(x_1,x_2)\|^2_{L^{-1}}
     + \ld  \la \lb M  -\ld \delta L \rb  \Delta {G}(x_1,x_2), \Delta {G}(x_1,x_2) \ra , \\
  \ea$ \\
   that is  (\ref{mohg}) with $i=2$\edem \\
%
\subsection{\label{djeu1} Proof of Proposition \ref{jeu1}.} According to  (\ref{subt1}), we have \\
 \ind $v_{n+1}=  (\ld w)  G_{\ld}^M(z_{n})= w \lb z_n -J_{\ld M^{-1} A} (z_n-\ld M^{-1} B(z_n)) \rb$, \\
 namely  \\
 \ind  $ (I+ \ld M^{-1} A)^{-1} (z_n - \ld  M^{-1}  B(z_n)) = z_n - \frac{1}{w} v_{n+1}$, \\
 which is equivalent to the inclusion \\
 \ind  $ M z_n - \ld  B(z_n) \in (M + \ld  A) \lb  z_n - \frac{1}{ w} v_{n+1}  \rb$. \\
 This, by $\ynp=  z_{n}-\frac{1}{w}v_{n+1} $ (according to (\ref{diss}))
  can be reduced to \\
   \ind $ M z_n - \ld B(z_n) \in  M z_n - \frac{1}{ w} M v_{n+1}+  \ld  A \lb   y_{n+1} \rb$, \\
  namely \\
  \ind   $ (\ld w)^{-1}  M v_{n+1} \in    (B(z_n) +  A (\ynp) )$, \\
  that is (\ref{piss0}).
  Hence, by \\
  \ind  $\ynp^*=  (\ld w)^{-1} M  v_{n+1} + B(\ynp)- B(z_n)  $  (according to (\ref{diss})), \\
  we get \\
  \ind $ \ynp^*   \in   B (\ynp) +  A (\ynp) $, that is  (\ref{piss1}). \\
   Next, by  definition of $\ynp^*$, we simply have
   \be \label{tuko}  \mbox{$\|\ynp^*\|_M\le   (\ld w)^{-1} \|M  v_{n+1} \|_M+   \|B(\ynp)- B(z_n)\|_M$}. \ee
   Let us estimate the two terms in the right side of the previous inequality. Concerning the first term, since  $M$ and $L$ are assumed to be bounded, by Remark \ref{pop} we have
   \be  \label{tuko1}  \mbox{$\|M  v_{n+1} \|^2_M=\la (M^{\frac{1}{2}})^6v_{n+1},v_{n+1}\ra =   \la M^{2} M^{1/2} v_{n+1}, M^{1/2} v_{n+1} \ra \le \|M\|^2  \| v_{n+1}\|_M^2$}. \ee
  Concerning the second term, we simply  get  \\
  \ind $\|B(\ynp)-B(z_n)\|_{M}\le \|M^{\frac{1}{2}}L^{\frac{1}{2}}\| \mul  \|B(\ynp)-B(z_n)\|_{L^{-1}},
 $ \\
 while  the co-coerciveness of $B$ (given by condition (\ref{cdi}))  yields   \\
  \ind $\|B(\ynp)-B(z_n)\|_{L^{-1}} \le \|L^{\frac{1}{2}}(\ynp-z_n)\|$,\\
  whence  it comes that
    \be  \label{tuko2}  \mbox{$\|B(\ynp)-B(z_n)\|_{M}\le \|M^{\frac{1}{2}}L^{\frac{1}{2}}\| \mul \|L^{\frac{1}{2}}(\ynp-z_n)\|$}, \ee
 where $ \ynp-  z_n=  -\frac{1}{ w} v_{n+1}$  (from (\ref{diss})).
   Then, by (\ref{tuko}), (\ref{tuko1}) and (\ref{tuko2}), we obtain
   \be \label{duke1}  \mbox{$\|\ynp^*\| \le   (\ld w)^{-1} \|M\|\mul  \|v_{n+1}\|_M + \frac{1}{ w} \|M^{\frac{1}{2}}L^{\frac{1}{2}}\|. \mul   \|L^{\frac{1}{2}}\vnp\|$}. \ee
   On the one hand, if   $M-\rho L$ is positive definite, we have  \\
   \ind   $\|\vnp\|_M^2 - \rho \|L^{\frac{1}{2}}\vnp\|^2 = \la \lb M- \rho L \rb \vnp, \vnp \ra \ge 0$, \\
   hence $  \|L^{\frac{1}{2}}\vnp\| \le  \rho ^{-1/2} \|\vnp\|_M$.
 On the other hand,  if   $M-\rho I$ is positive definite, we get   \\
 \ind $\|M^{\frac{1}{2}}\vnp\|^2 - \rho \|\vnp\|^2 = \la \lb M- \rho I \rb \vnp, \vnp \ra \ge 0$, \\
 which yields  $  \|\vnp\| \le  \rho ^{-1/2} \|\vnp\|_M$, hence \\
 \ind $\|L^{\frac{1}{2}} \vnp \|\le  \rho ^{-1/2} \|L^{\frac{1}{2}}\| .  \|\vnp\|_M$. \\
 Consequently, regarding the previous two situations, by (\ref{duke1})  we are led to \\
  \ind $\|\ynp^*\|_M \le   (\ld w)^{-1} \|M\|\mul  \|v_{n+1}\|_M + \frac{1}{ w} \|M^{\frac{1}{2}}L^{\frac{1}{2}}\|\mul
    \rho ^{-1/2} (1+\|L^{\frac{1}{2}}\|) \mul  \|\vnp\|_M$, \\
    which  amounts to (\ref{piss2})\edem \\
\subsection{\label{pap2} Proof of Proposition \ref{gfru0}.} For the sake of simplicity, we write $\la ., . \ra $ instead of $\la ., . \ra _M$.
 Setting $\tau_n:=e+ \nu_{n+1}$, we observe  that    (\ref{fdc})  can be alternatively  expressed as
   \be
    \label{hypg1}
  \mbox{$\ttn =\frac{\nu_n }{\tau_n }$}. \ee
 As another crucial parameter arising in our study, we consider   the real sequence  $(\gamos)$ (with $s > 0$) defined by
  \be \label{tld} \mbox{$\gamos=1- \frac{s}{\tau_n}$}. \ee
 The following  elementary observation will  be particularly helpful for  the sequel of our study. \\
  \begin{remark} \label{sld} {
 For $s \in ( 0, e ]$, we have $(\gamos) \subset (0,1)$. \\ \\
} \end{remark}
 It is  readily  noticed   that  $\dot{F}_{n}(s,q)$ can be formulated as
\be \label{est1} \ba
   \ind \mbox{$\dot{F}_{n+1}(s,q) =  s  (\nu_{n+1} {a}_{n+1}- \nu_{n}  a_n)  +  (s e) \dot{b}_{n+1} + \nu^2_{n+1}{c}_{n+1}-  \nu_{n}^2{c}_{n} $},
   \ea \ee
where $a_n:=\la x_n-q,   \dxn \ra$,   $b_n := (1/2) \|x_n-q\|^2$ and  $c_n :=(1/2) \|  \dxn\|^2$. \\
   Note also that for any bilinear symmetric form $\la .,.\ra_E$ on a real vector space $E$ and for any sequences
$\{  \phi_{n} , \varphi_{n} \} \subset E$ we have the discrete derivative rules:
\begin{subeqnarray}
\label{ddr}
\label{ddr1}
&& \la \phi_{n+1} , \varphi_{n+1} \ra _E- \la \phi_{n} , \varphi_{n} \ra_E=\la  \dot{\phi}_{n+1}, \varphi_{n+1} \ra_E + \la  \phi_{n} ,  \dot{\varphi}_{n+1} \ra_E,  \\
\label{ddr2}
&& \la \phi_{n+1} , \varphi_{n+1} \ra_E - \la \phi_{n} , \varphi_{n} \ra_E=\la  \dot{\phi}_{n+1}, \varphi_{n} \ra_E + \la  \phi_{n+1} ,  \dot{\varphi}_{n+1} \ra_E.
\end{subeqnarray}
The sequel of the proof can be divided into the following parts {\bf (1)-(5)}: \\
 \\
\ind {\bf (1)} {\bf Basic estimates.} Setting     $\tau_n=e+ \nu_{n+1}$,      $P_n= \la q-\xnp, \dxnp  \ra$ and  $ R_n=  \la q-\xnp  , \dxn \ra$,  we  establish  the elementary but useful facts  below:
\begin{subeqnarray}
 \label{cgen}
 \label{pmr1}
   \dot{a}_{n+1}  & =&  \la  \dxnp  , \dxn  \ra  -  P_n +  R_n, \\
\label{pmr2}
    \dot{b}_{n+1} & =&-P_n  -
  (1/2)\| \dxnp\|^2, \\
   \label{cgen2} \nu_{n+1}
    a_{n+1} -\nu_n a_n  &=&   \nu_n   \la \dxnp, \dxn  \ra -    \nu_{n+1} P_n +  \tau_n \lb     P_n  -   \la d_{n}, \xnp-q \ra \rb .
\end{subeqnarray}
    \ind  Let us prove (\ref{pmr1})- (\ref{pmr2}). From $a_n:=\la x_n-q,   \dxn \ra$, by the  rule (\ref{ddr2})
we simply  have
\[ \mbox{$\dot{a}_{n+1}  = \la  \dxnp  , \dxn \ra + \la \xnp -q ,  \dxnp -\dxn \ra$}, \]
which gives us  (\ref{pmr1}).
 From  $ 2  b_{n+1}= \|\xnp-q\|^2$, by the derivative  rule (\ref{ddr1}) we get
 \[ \mbox{$ \ba 2 \dot{b}_{n+1}   =\la \dxnp, \xnp-q \ra + \la \xn -q, \dxnp \ra \\
 \hskip 1.cm =\la \dxnp, \xnp-q \ra + \la x_n -\xnp , \dxnp \ra +  \la \xnp -q, \dxnp \ra, \ea $} \]
     namely $ 2 \dot{b}_{n+1}     = -2 P_n  - \| \dxnp\|^2$,
which leads to  (\ref{pmr2}). \\
\\
\ind Let us prove (\ref{cgen2}). From   $a_n= \la q-x_n, -   \dxn \ra $, we  simply get \\
\ind $\ba a_n=\la q-\xnp,     -\dxn \ra + \la \dxnp, -\dxn  \ra=- R_n  - \la \dxnp,    \dxn\ra,
    \ea$ \\
    while the   derivative  rule (\ref{ddr2}) yields
    \ind $ \nu_{n+1} a_{n+1} -\nu_n a_n= \dot{\nu}_{n+1}  a_{n} + {\nu}_{n+1}\dot{a}_{n+1}$. \\
This, in light of   (\ref{pmr1}), amounts to
  \be \label{rat}  \mbox{$\ba  a_{n+1}\nu_{n+1} - a_n \nu_n  = \dot{\nu}_{n+1} \lb - R_n  - \la \dxnp,    \dxn\ra\rb +  {\nu}_{n+1}(  \la  \dxnp  , \dxn  \ra  -  P_n +  R_n) \\
\hskip 2.5cm  = \nu_n \la \dxnp, \dxn  \ra -   \nu_{n+1} P_n +    \nu_{n}  R_n.\ea $} \ee
  Furthermore,     (\ref{adr0})  gives us   $ \dxnp  +   d_{n} -  \theta_n  \dxn =0$.
Taking the scalar product of each side of this equality by $q-\xnp$ yields \\
 \ind $\ttn  R_n   =      P_n  -     \la d_{n}, \xnp-q \ra$. \\
So, noticing that   $ \nu_{n}  =\tau_n \ttn$, we get \\
 \ind $  \nu_{n} R_n =\tau_n  (\ttn R_n)=\tau_n \lb     P_n  - \la d_{n}, \xnp-q \ra    \rb$, \\
which, in light of (\ref{rat}),    entails
  (\ref{cgen2}). \\
\\
\ind {\bf (2)} {\bf An estimate from the inertial part.}
 Now, given  $(s,q) \in [0,\li) \mul \Hs $, we    prove  that the discrete  derivative $  \dot{F}_{n+1}(s,q)$ satisfies
 \be  \label{fgen}
  \hskip .5cm  \ba \dot{F}_{n+1}(s,q)
      +   s (e+ \nu_{n+1}) \la d_{n}, \xnp -q  \ra  \\
      \hskip 1.cm   =  s \nu_n  \la \dxnp, \dxn \ra
     - \frac{1}{2} \lb s e  -  \nu_{n+1}^2    \rb  \|\dxnp\|^2 -  \frac{1}{2}   \nu_n ^2  \| \dxn  \|^2.
    \ea   \ee
 \ind  Indeed, in light of  (\ref{est1}) and  (\ref{cgen}), we obtain
 \[ \mbox{ $ \ba \dot{F}_{n+1}  = s \lb  \nu_n   \la \dxnp, \dxn  \ra -    \nu_{n+1} P_n +  \tau_n \lb     P_n  -    \upsilon_n \rb  \rb  \\
 \ind \ind \ind + (se) \lb -P_n  -
   \frac{1}{2}\| \dxnp\|^2 \rb +  \frac{1}{2} \lb \nu_{n+1}^2   \| \dxnp \|^2-
     \nu_{n}^2  \|\dxn\|^2 \rb \\
   \ind \ind =   s   \nu_n   \la \dxnp, \dxn  \ra + s(-  \nu_{n+1}+  \tau_n-e) P_n +  \frac{1}{2} \lb se - \nu_{n+1}^2   \rb \|\dxnp\|^2 \\
  \hskip 1.cm -s \tau_n  \la d_{n}, \xnp-q \ra ,
   \ea$} \]
 hence, by  rearranging the terms in the previous equality, we are  led to
   \[ \mbox{ $ \ba \dot{F}_{n+1} + s  \tau_n  \la d_{n}, \xnp -q  \ra  =  s \nu_n  \la \dxnp, \dxn \ra
     + \frac{1}{2} \lb \nu_{n+1}^2 - s e        \rb  \|\dxnp\|^2 -  \frac{1}{2}   \nu_n ^2  \| \dxn  \|^2.   \ea$} \]
    This obviously leads to the desired equality. \\
    \\
\ind {\bf (3)} {\bf An estimate from the proximal part.} We prove   that,  for any $\gamc_n \neq 1$,  it holds that
 \be \label{sup2i} \ba   \gamc_n    \la  d_n,\dxnp \ra   +    \frac{ 1 }{2}   \| \dxnp- \ttn \dot{x}_n \|^2 \\
 \hskip 2.cm
   =    - \ttn \bgamc \la \dot{x}_n,\dxnp\ra +  \frac{1}{2}    \ttn^2 \|\dot{x}_n \|^2   -  \lb \gamc_n- \frac{1}{2}  \rb \|\dxnp\|^2.
    \ea \ee
  \ind Indeed, we have   $\dxnp=\ttn \dot{x}_n -d_n$ (from (\ref{adr0})), hence, for any  $\gamc_n \neq 1$, and setting $H_n=  \dxnp - \bgamc^{-1}\ttn \dot{x}_n $, we get \\
  \ind $ \bgamc H_n=  \bgamc\dxnp - \ttn \dot{x}_n  = -d_n- \xi_n \dxnp$, \\
  or equivalently
  \be  \label{stor1} \ba  \ind \ind \gamc_n   \dxnp     = -\bgamc H_n  -   d_n . \ea \ee
   Furthermore, by $-d_n= \dxnp- \ttn \dot{x}_n$ (again using  (\ref{adr0})), we simply obtain
   \be \label{stor2} \ba
  \la (-d_n), H_n \ra  =    \la  \dxnp- \ttn \dxn,  \dxnp  - \bgamc^{-1}\ttn \dxn \ra \nonumber  \\
    \hskip 2.cm  =    \|\dxnp\|^2 + \bgamc^{-1}\ttn^2 \|\dxn\|^2 - \mbox{$ \frac{2-\xi_n}{\bgamc}\ttn    \la \dxnp, \dxn\ra$.}  \
   \ea \ee
  Therefore,   taking  the scalar product of each  side of  (\ref{stor1}) with $d_n$,  and  adding  $(1/2) \| d_n\|^2$ to the resulting equality, next  using (\ref{stor2})  and   $ \| d_n\|^2= \|\dxnp- \ttn \dxn \|^2$  we get
  \[ \ba    \gamc_n  \la d_n , \dxnp  \ra + \frac{1}{2} \| d_n\|^2  = \bgamc \la (-d_n), H_n \ra   -  \frac{1}{2} \|  d_n\|^2
 \\
  \hskip 3.cm = \bgamc \lb \|\dxnp\|^2 + \frac{ \ttn^2 }{\bgamc} \|\dxn \|^2 - \frac{2-\xi_n}{\bgamc} \ttn \la \dxnp, \dxn\ra  \rb  \\
  \hskip 5.cm-  \frac{1}{2}\lb \|\dxnp\|^2 + \ttn^2 \| \dxn \|^2 - 2 \ttn \la \dxnp, \dxn\ra\rb \\
    \hskip 3.cm =  -(1-\xi_n) \ttn \la \dxnp, \dxn\ra + \frac{1}{2}\ttn^2 \| \dxn \|^2 + \lb \frac{1}{2} -\gamc_n \rb  \|\dxnp\|^2 .
 \ea \]
 This  yields (\ref{sup2i}). \\
 \\
 \ind {\bf (4)} {\bf Combining proximal and  inertial effects}. At once we show  for $s \in \lb 0, e \right ]$ that  the iterates verify (for $n \ge 0$)
  \be \label{fed} \ba
 \dot{F}_{n+1}(s,q)  + \frac{1}{2}  \tau_n^2 \|W_n\|^2 \\
  + (s \tau_n)\la d_n, \xnp-q \ra   +  \gamos \tau_n^2 \la  d_n,\dxnp  \ra
       =  -T_n(\dxnp) ,
    \ea \ee
    where $\tau_n=e+\nu_{n+1}$, $\gamos=1- s \tau_n ^{-1}$, $W_n=\dxnp - \ttn \dxn $, while   $T_n(x)$ is defined for any
    $x \in \Hs$ by
 \be \label{dtn} \mbox{$T_n(x) =
  \frac{1}{2} \lb se  - \nu_{n+1}^2
  +   \tau_n^2  \lb 2 \gamos - 1  \rb
 \rb   \|x\|^2$}. \ee
    \ind Indeed, from  (\ref{fgen}) we know that
 \be  \label{floqas}
   \ba \dot{F}_{n+1}(s,q)    + (s \tau_n)\la d_n, \xnp-q \ra   \\
   \hskip .5cm   =   - \lb - s \nu_n  \la \dxnp, \dxn \ra
    +  \frac{1}{2}   \nu_n ^2  \| \dxn  \|^2 + \frac{1}{2} \lb s e  -  \nu_{n+1}^2    \rb  \|\dxnp\|^2  \rb .
    \ea   \ee
    Furthermore, by   $\gamos=1- s \tau_n^{-1} $ and assuming that $s \in \lb 0, e \right ]$,   we have  $(\gamos) \subset (0,1)$ (according to   Remark \ref{sld}).  So, from  (\ref{sup2i}), we also have
      \be \label{nsup2i} \ba    \gamos   \la d_n ,\dxnp \ra   +  \frac{ 1 }{2}  \|W_n\|^2  \\
      \hskip 1.cm   = -\lb  \ttn (1-\gamos) \la \dxn, \dxnp \ra- \frac{1}{2}  \ttn^2 \|\dxn\|^2   +   \lb \gamos -\frac{1}{2} \rb \|\dxnp\|^2 \rb  .  \ea\ee
     Then  multiplying equality (\ref{nsup2i})  by $\tau_n^2$ (while  noticing that $  \tau_n^2=(s\tau_n) (1-\gamos)^{-1})$,
       and adding the resulting equality to (\ref{floqas})  gives us
 \be  \label{floqas2}  \ba  \gamos \tau_n^2    \la d_n ,\dxnp \ra    + \frac{1}{2} \tau_n^2 \| W_n \|^2 \\
  +  \dot{F}_{n+1}(s,q)   + (s \tau_n)\la d_n , \xnp-q \ra
     = -\lb w_{n}   \la  \dxn, \dxnp \ra
     +\eta_{n} \| \dxn \|^2+ \sigma_{n}  \|\dxnp\|^2 \rb, \ea  \ee
    together with  the following parameters
\[ \ba  \mbox{$ {w}_{n} = - s \nu_n +  \tau_n^2  \ttn (1-\gamos)  $}, 
\ind \ind \mbox{$  {\eta}_{n}=\frac{1}{2}   \nu_n ^2 -  \frac{1 }{2 }  \tau_n^2  \theta_n^2$},  \ind 
\ind  \mbox{$ {\sigma}_{n} =\frac{1}{2} \lb s e  -  \nu_{n+1}^2    \rb  + \tau_n^2    \lb  \gamos-  \frac{1}{2}   \rb$.} \ea \]
 Clearly, by   $\tau_n \ttn =\nu_n $ (from (\ref{fdc})),  and noticing that $(1-\gamos)^{-1}=\frac{\tau_n}{s}$,  we are led to
 \[ \ba  \ind \mbox{$  {w}_{n} =  - s \nu_n
 +  s \nu_n   =  0$}, \\
 \ind \mbox{$  {\eta}_{n}= \frac{1}{2}  \nu_n ^2    -  \frac{  \tau_n^2 \theta_n^2}{2 }= \frac{1}{2} \lb \nu_n^2  -   \nu_n^2 \rb=0$},  \\
 \ind \mbox{${\sigma}_{n} = \frac{1}{2} \lb se
  -  \nu_{n+1}^2  \rb
    +    \tau_n^2\lb \gamos - \frac{1}{2}  \rb = \frac{1}{2} \lb se  -\nu_{n+1}^2
  +  \tau_n^2  \lb 2 \gamos - 1  \rb
 \rb $}. \ea  \]
 This leads to (\ref{fed})-(\ref{dtn}). \\
 \\
\ind {\bf (5)} It remains to   simplify the  formulation  given by  (\ref{dtn}).
   A simple computation yields \\
 \ind   $2 \gamos -1= 1-  \frac{2s}{(e+ \nu_{n+1})}$ \es \es (as  $\gamos=1- \frac{s}{e+ \nu_{n+1}}$), \\
   which implies that
  \[ \mbox{$ \ba ( e+ \nu_{n+1}) ^2 (2 \gamos -1) = \lb e^2+ 2  e \nu_{n+1}  + ( \nu_{n+1})^2 \rb-  2s  \lb e+  \nu_{n+1}   \rb \\
  \hskip 3.cm  = e   \lb  e -s \rb - s e + 2  \nu_{n+1}  \lb e -s \rb
  + ( \nu_{n+1} )^2  \\
   \hskip 3.cm =  \lb e + 2  \nu_{n+1}\rb  \lb e -s \rb - s e +  (  \nu_{n+1})^2.
   \ea$} \]
  As a consequence, by   (\ref{dtn}) we obtain \\
   \ind $T_{n}(x)= \frac{\lb e  -s \rb}{2}  \lb e + 2  \nu_{n+1}\rb    \|x\|^2$. \\
    This ends the proof\edem   
%
%
 
\end{document}